\title{On a hypergraph Tur\'an problem of Frankl
}
\author{
Peter Keevash \thanks{
Princeton University, Princeton, NJ 08540, USA. Email address:
keevash@math.princeton.edu } \and
Benny Sudakov \thanks{Department of Mathematics,
Princeton University, Princeton, NJ 08540, USA.
Email address: bsudakov@math.princeton.edu.
Research supported in part by NSF grant DMS-0106589.}
        }
\date{}
\documentclass [11pt]{article}
\usepackage{amsfonts}
\usepackage{amsmath}
\newtheorem{theo}{Theorem}[section]
\newtheorem{prop}[theo]{Proposition}
\newtheorem{lemma}[theo]{Lemma}
\newtheorem{coro}[theo]{Corollary}

\newtheorem{claim}[theo]{Claim}

\begin{document}
\maketitle

\begin{abstract}
Let ${\cal C}^{(2k)}_r$ be the $2k$-uniform hypergraph obtained by
letting $P_1,\cdots,P_r$ be
pairwise disjoint sets of size $k$ and taking as edges all sets
$P_i \cup P_j$ with $i \neq j$. This can be thought of as the
`$k$-expansion' of the complete graph $K_r$: each vertex has
been replaced with a set of size $k$. An example of a
hypergraph with vertex set $V$ that does not contain ${\cal C}^{(2k)}_3$ can be 
obtained by partitioning $V = V_1 \cup V_2$ and taking as edges all sets
of size $2k$ that intersect each of $V_1$ and $V_2$ in an odd number
of elements. Let ${\cal B}^{(2k)}_n$ denote a hypergraph on $n$ vertices
obtained by this construction that has as many edges as possible.
We prove a conjecture of Frankl, which states that any hypergraph
on $n$ vertices that contains no ${\cal C}^{(2k)}_3$ has at most as
many edges as ${\cal B}^{(2k)}_n$.

Sidorenko has given an upper bound of $\frac{r-2}{r-1}$ for
the Tur\'an density of ${\cal C}^{(2k)}_r$ for any $r$, and a construction
establishing a matching lower bound when $r$ is of the form
$2^p+1$. In this paper we also show that when $r=2^p+1$, any ${\cal 
C}^{(4)}_r$-free hypergraph of density $\frac{r-2}{r-1}-o(1)$ looks 
approximately like
Sidorenko's construction. On the other hand, 
when $r$ is not of this form, we show that corresponding constructions do not 
exist and 
improve the upper bound on the Tur\'an density of ${\cal C}^{(4)}_r$ to
$\frac{r-2}{r-1} - c(r)$, where $c(r)$ is a constant
depending only on $r$.

The backbone of our arguments is a strategy
of first proving approximate
structure theorems, and then showing that any imperfections in the
structure must lead to a suboptimal configuration. The tools for
its realisation draw on extremal graph theory, linear algebra, the 
Kruskal-Katona theorem
and properties of Krawtchouck polynomials.
\end{abstract}

\section{Introduction}
Given an $r$-uniform hypergraph $\cal F$, the Tur\'an number $ex(n, {\cal F})$ 
of $\cal F$ is
the maximum number of edges in an $r$-uniform hypergraph on $n$ vertices 
that does not contain a copy of $\cal F$. Determining these numbers is one of 
the main
challenges in Extremal Combinatorics. For ordinary graphs (the case $r=2$) a 
rich
theory has been developed, initiated by Tur\'an in 1941, who 
solved the problem for complete graphs. He also posed the question of 
finding $ex(n, {\cal K}_s^{(r)})$ for 
complete hypergraphs with $s>r>2$, but to this day not one single instance 
of this problem has been solved.
It seems hard even to determine the {\em Tur\'an density}, which for general
$\cal F$ is defined as $\pi({\cal F}) = \lim_{n \rightarrow \infty} ex(n,{\cal 
F})/{n \choose r}$. 
The problem of finding the numbers $ex(n, {\cal F})$ when $r>2$ is notoriously 
difficult, and 
exact results on hypergraph Tur\'an numbers are very rare (see
\cite{F1,Si1} for surveys).
In this paper we obtain such a result for a sequence of hypergraphs introduced 
by Frankl.

Let ${\cal C}^{(2k)}_r$ be the $2k$-uniform hypergraph obtained by
letting $P_1,\cdots,P_r$ be
pairwise disjoint sets of size $k$ and taking as edges all sets
$P_i \cup P_j$ with $i \neq j$. This can be thought of as the
`$k$-expansion' of the complete graph $K_r$: each vertex has
been replaced with a set of size $k$. The Tur\'an problem
for ${\cal C}^{(2k)}_3$ was first considered by Frankl \cite{Fr}, who
determined the density $\pi({\cal C}^{(2k)}_3)=1/2$.

Frankl obtained a large
${\cal C}^{(2k)}_3$-free hypergraph on $n$ vertices by partitioning an 
$n$-element set $V$ into $2$ parts $V_1,V_2$ and taking those edges which
intersect each part $V_i$ in an odd number of elements. When
the parts have sizes $\frac{n}{2} \pm t$ we denote this hypergraph
by ${\cal B}^{(2k)}(n,t)$. To see that it is ${\cal C}^{(2k)}_3$-free, consider
any $P_1,P_2,P_3$ that are pairwise disjoint sets of $k$ vertices.
Then $|V_1 \cap P_i|$ and $|V_1 \cap P_j|$ have the same parity for
some pair $ij$, so $P_i\cup P_j$ is not an edge. 
Let $t^*$ be chosen to maximise the number of edges in ${\cal B}^{(2k)}(n,t)$,
and denote any hypergraph obtained in this manner by ${\cal B}^{(2k)}_n$.
Write $b_{2k}(n)$ for the number of edges in ${\cal B}^{(2k)}_n$. Frankl 
\cite{Fr} conjectured that the maximum number of edges in a ${\cal 
C}^{(2k)}_3$-free hypergraph is always achieved by some ${\cal B}^{(2k)}_n$. 
Our first theorem proves this conjecture.

\begin{theo}
\label{frankl}
Let $H$ be a $2k$-uniform hypergraph on $n$ vertices that does not contain a
copy of ${\cal C}^{(2k)}_3$ and let $n$ be sufficiently large. Then the number 
of
edges in $H$ is at most $b_{2k}(n)$, with equality only when $H$ is
a hypergraph of the form ${\cal B}^{(2k)}_n$.
\end{theo}

The proof of this theorem falls naturally into two parts. The first
stage is to prove a `stability' version, which is that any hypergraph
with close to the maximum number of edges looks approximately like
some ${\cal B}^{(2k)}(n,t)$. Armed with this, we can analyse any imperfections
in the structure and show that they must lead to a suboptimal
configuration, so that the optimum is indeed achieved by the
construction. This strategy was also used recently in \cite{KS} to prove the
conjecture of S\'os on the Tur\'an number of the Fano plane,
so this seems to be a useful tool for developing the Tur\'an theory
of hypergraphs.

For general $r$, Sidorenko \cite{Si} showed that the Tur\'an density of ${\cal 
C}^{(2k)}_r$ 
is at most $\frac{r-2}{r-1}$. This is a consequence of Tur\'an's theorem 
applied
to an auxiliary graph $G$ constructed from a $2k$-uniform hypergraph $H$; 
the vertices of $G$ are the $k$-tuples of vertices of $H$, and
two $k$-tuples $P_1$,$P_2$ are adjacent if $P_1 \cup P_2$ is an edge of $H$.
He also gave a construction for a matching lower bound when
$r$ is of the form $2^p+1$, which we now describe.
Let $W$ be a vector space of dimension $p$ over the field 
$GF(2)$, i.e. the finite field with $2$
elements $\{0,1\}$. Partition a set of vertices $V$ as $\bigcup_{w \in W} 
V_w$. Given $t$ and a $t$-tuple of vertices $X = x_1 \cdots x_t$ with $x_i \in 
V_{w_i}$ 
let $\Sigma X = \sum_1^t w_i$. Define a $2k$-uniform hypergraph $H$, where
a $2k$-tuple $X$ is an edge iff $\Sigma X \neq 0$. Observe that this
doesn't contain a copy of ${\cal C}^{(2k)}_r$. Indeed, if 
$P_1,\cdots,P_r$ are disjoint $k$-tuples
then there is some $i \neq j$ with $\Sigma P_i = \Sigma P_j$ (by the
pigeonhole principle). Then $\Sigma (P_i \cup P_j) =\Sigma P_i+ \Sigma P_j
= 0$, so $P_i \cup P_j$ is
not an edge. To see that this construction can achieve the stated Tur\'an 
density,
choose the partition so that $|V_w| =|V|/(r-1)$.
Then a random (average) $2k$-tuple is an edge with probability 
$\frac{r-2}{r-1} + o(1)$, as
can be seen by conditioning on the positions of all but one element.

This construction depends essentially on an algebraic structure, which
only exists for certain values of $r$. We will
show that this is an intrinsic feature of the problem, by proving
a stronger upper bound on the Tur\'an density of ${\cal C}^{(4)}_r$ when
$r$ is not of the form $2^p+1$. 

\begin{theo}
\label{c4r}
Suppose $r \geq 3$, and let $H$ be a $4$-uniform hypergraph on $n$
vertices with at least $\big(\frac{r-2}{r-1} - 10^{-33}r^{-70} \big) {n
\choose 4}$ edges. If $H$ is ${\cal C}^{(4)}_r$-free, then $r=2^p+1$ for
some integer $p$.
\end{theo}

In contrast to Theorem \ref{frankl} this is a result showing that certain
constructions do not exist, so it is perhaps surprising
that its proof also uses a stability argument. We study the
properties of a ${\cal C}^{(4)}_r$-free hypergraph with density close
to $\frac{r-2}{r-1}$ and show that it give rise to the edge coloring
of the complete graph $K_{r-1}$ with special properties. Next we prove that 
for such edge-coloring there is a natural $GF(2)$ vector space
structure on the colors. Of 
course, such a space has cardinality $2^p$, for some $p$, so we get a
contradiction unless $r=2^p+1$.

A complication arising in Theorem \ref{frankl} is that the optimum construction
is not achieved by a partition into two equal parts. Finding $t$ to maximise
the number of edges in ${\cal B}^{(2k)}(n,t)$ 
is an interesting problem in enumerative combinatorics, equivalent
to finding the minima of binary Krawtchouk polynomials. This is a family
of polynomials orthogonal with respect to the uniform
measure on a $n$-dimensional cube that play an important r\^ole in the analysis
of binary Hamming association schemes (see, e.g., \cite{KL}). Despite some 
uncertainty in the location of their minima, the known bounds are sufficient 
for us to show that some ${\cal B}^{(2k)}(n,t)$ must be optimal.

In the case $k=2$ one can compute the size of ${\cal B}^{(2k)}(n)$ precisely, 
and there are 
considerable
simplifications of the argument, so in the next section for illustrative
purposes we start by giving a separate proof for this case.
Section 3 contains a stability theorem for ${\cal C}^{(2k)}_3$ and 
the general case of Theorem \ref{frankl}. Then in Section 4 we prove a 
stability result for ${\cal C}^{(4)}_r$ for all $r$, and use it to establish 
Theorem \ref{c4r}.
The final section of the paper contains some concluding remarks.

We will assume throughout this paper that $n$ is sufficiently large.

\section{The Tur\'an number of ${\cal C}^{(4)}_3$}

We start by proving Frankl's conjecture for $4$-uniform hypergraphs. This 
will serve to illustrate our method, as it has fewer complications than 
the general case. In addition, in this case it is easy to compute the 
Tur\'an numbers of ${\cal C}^{(4)}_3$ precisely.

We recall that ${\cal C}^{(4)}_3$ is the $4$-uniform hypergraph 
with three edges $\{abcd,abef,cdef\}$. We can obtain a large
${\cal C}^{(4)}_3$-free graph on $n$ vertices by partitioning an 
$n$-element
set into $2$ parts and taking those edges
which have $1$ point in either class and $3$ points in the other.
To see this, think of an edge as being the union of $2$
different {\em types} of pairs of 
vertices: one type consisting of pairs with 
both vertices in one class, the other consisting
of pairs that have one point of each class.
Given any $3$ pairs there are $2$ of the same type, and these do not
form an edge in the construction.

To maximise the number of edges in this bipartite construction, it is 
{\em not} the case
that the two parts have sizes as equal as possible, but we will see that
the difference in the sizes should be at most of order $\sqrt{n}$.
Let ${\cal B}(n,t)$ denote the $4$-uniform hypergraph obtained 
by partitioning an $n$-element set into $2$ parts
with sizes $\frac{n}{2}+t$ and $\frac{n}{2}-t$,
and taking those edges
which have $1$ point in either class and $3$ points in the other.
Let $b(n,t)$ be the number
of edges in ${\cal B}(n,t)$ and let $d(n,t)$ be the degree of any vertex
belonging to the side with size $\frac{n}{2}+t$. Then the
vertices on the side with size $\frac{n}{2}-t$ have degree
$d(n,-t)$. We will start with some estimates on these parameters.
By definition, 
\begin{eqnarray}
\label{heq}
b(n,t) & = & \Big(\frac{n}{2}+t\Big){{\frac{n}{2}-t} \choose 3} +
             \Big(\frac{n}{2}-t\Big){{\frac{n}{2}+t} \choose 3} 
        =  \frac{n^4-6n^3+8n^2-16t^4-32t^2+24t^2n}{48} \nonumber\\
       & = & \frac{1}{48}\Big( \Big(n^2-3n+4\Big)^2 -
    \Big(4t^2-3n+4 \Big)^2 \Big)\,.
\end{eqnarray}
Thus to maximise $b(n,t)$ we should pick a value of $t$
that minimises $4t^2-3n+4$, subject to the restriction that
when $n$ is even $t$ has to be an integer, and
when $n$ is odd $t+\frac{1}{2}$ has to be an integer.
Let ${\cal B}_n$ denote a hypergraph ${\cal B}(n,t^*)$, where $t^*$ 
is such a value of $t$. By symmetry we can take $t^*>0$. There is
usually a unique best choice of $t^*$, but for
some $n$ there are $2$ equal choices of $t^*$. Note
that for any best choice we certainly have 
$\big|t^*-\sqrt{3n/4-1}\big|\leq 1/2$.

Let $b(n)$ be the number of edges in ${\cal B}_n$. Then
$$\big|48b(n) - (n^2-3n+4)^2\big| = \big|4(t^*)^2 - 3n + 4\big|^2
< 50n\,.$$
It will be useful later to consider the following estimate
which  follows immediately from the last
inequality for sufficiently large $n$
\begin{equation}
\label{diff}
b(n)-b(n-1) > \frac{1}{12}n^3 - \frac{1}{2}n^2\,.
\end{equation}
Next we give an explicit formula for the degrees in ${\cal B}(n,t)$
\begin{equation}
\label{degree}
d(n,t)  =  \left(\frac{n}{2} - t\right){{n/2 + t - 1} \choose 2} +
 {{n/2 - t} \choose 3} 
 =  \frac{n^3-6n^2+8n+12t^2}{12} + \frac{6tn-8t^3-16t}{12}\,.
\end{equation}
We finish these calculations with an upper bound on the
 maximum degree of ${\cal B}_n$
\begin{equation}
\label{dmax}
\Delta(n)  =  \frac{1}{12}\big(n^3-6n^2+8n+12(t^*)^2\big) + 
\frac{1}{12}\big|6t^*n-8(t^*)^3-16t^*\big|
< \frac{1}{12}n^3 - \frac{1}{2}n^2 + n^{3/2}\,.
\end{equation}

The first step in the proof is to show that any 
${\cal C}^{(4)}_3$-free
$4$-uniform hypergraph $H$ with density close to $1/2$ has the correct
approximate structure. To do so we need a few definitions. If we have a 
partition of the vertex set of $H$
as $V(H)=V_1 \cup V_2$ we call a $4$-tuple of vertices {\em good} if
it has either $1$ point in $V_1$ and $3$ points in $V_2$ or $1$ point in $V_2$
and $3$ points in $V_1$; otherwise we call it {\em bad}. With respect to
$H$, we call a $4$-tuple {\em correct} if it is either a good edge or
a bad non-edge; otherwise we call it {\em incorrect}.
We obtain the following stability result.
\begin{theo} 
\label{stability}
For every $\epsilon>0$ there is $\eta>0$ so that if $H$ is a 
${\cal C}^{(4)}_3$-free
$4$-uniform hypergraph with $e(H)>b(n)-\eta n^4$ then there
is a partition of the vertex set as $V(H)=V_1 \cup V_2$ such that all but
$\epsilon n^4$ $4$-tuples are correct.
\end{theo}

In the proof of this result we need a special case of
the Simonovits stability theorem \cite{S1} for graphs, which we 
recall. It states that
for every $\epsilon'>0$ there is $\eta'>0$ such that if $G$ is a
triangle free graph on $N$ vertices with at least $(1-\eta'){N \choose 
2}/2$ edges then there is a partition of the vertex set as
$V(G) = U_1 \cup U_2$ with $e_G(U_1)+e_G(U_2)<\epsilon' N^2$. 
 
\vspace{0.1cm}
\noindent 
{\bf Proof of Theorem \ref{stability}.}\
Define an auxiliary graph $G$ whose vertices are all pairs of vertices of 
$H$, and where the pairs $ab$ and $cd$ are adjacent exactly when $abcd$ is 
an edge of $H$. Since $H$ is ${\cal C}^{(4)}_3$-free we see that $G$ is 
triangle-free.
Also, each edge of $H$ creates exactly $3$ edges in $G$ (corresponding
to the $3$ ways of breaking a $4$-tuple into pairs) so
$$e(G) > 3\Big(b(n)-\eta n^4\Big) > \big(1-50\eta\big) \frac{1}{2} {{n 
\choose 2} \choose 2}\,.$$

Choose $\eta$ so that Simonovits stability applies with $\eta'=50\eta$,
$N={n \choose 2}$ and $\epsilon'=\epsilon^2/500$. We can
also require that $\eta < \epsilon^2/500$. We get a partition of the
pairs of vertices of $H$ as $U_1 \cup U_2$, where all but
$\epsilon'N^2<\epsilon^2 n^4/2000$ edges of $H$ are formed by taking
a pair from $U_1$ and a pair from $U_2$. 

We will call the pairs in $U_1$ {\em red}, and the pairs in $U_2$ {\em 
blue}. A $4$-tuple $abcd$ will be called {\em properly coloured} if 
either 

\noindent
(i) $abcd$ is an edge of $H$ and each of the $3$ sets 
$\{ab,cd\}$,$\{ac,bd\}$,$\{ad,bc\}$
has one red pair and one blue pair, or

\noindent
(ii) $abcd$ is not an edge and each of the $3$ sets 
$\{ab,cd\}$,$\{ac,bd\}$,$\{ad,bc\}$
consists of two pairs with the same colour.

An improperly coloured $4$-tuple is either an edge that is the union of two
pairs of the same colour or
a non-edge which is the union of two pairs with different colours.
There are at most $\epsilon^2 n^4/2000$ of the former $4$-tuples,
and the number of latter is at most
$$|U_1||U_2| - \Big(e(G)-\epsilon'N^2 \Big)\leq 
\frac{50\eta}{2}\frac{N^2}{2} +\epsilon'N^2 \leq
\left(\frac{50\eta}{16}+\epsilon'/4\right)n^4<\epsilon^2 n^4/140\,.$$
 Therefore all but 
$\big(\epsilon^2/140+\epsilon^2/2000\big)n^4<\epsilon^2n^4/130$ $4$-tuples 
are properly coloured.

A simple counting argument shows that there is a pair $ab$
so that for all but 
${4 \choose 2}\big(\epsilon^2n^4/130\big)/{n \choose 2}<\epsilon^2 
n^2/10$ other 
pairs $cd$ the
$4$-tuple $abcd$ is properly coloured. Without loss of generality $ab$ is red.
Partition the vertices of $V-ab$ into
$4$ sets according to the colour of the edges they send to $\{a,b\}$. We label
these sets $RR$,$BB$,$RB$,$BR$, where $R$ means `red', $B$ means `blue' and a
vertex $c$ belongs to the set that labels the colours of the edges $ca,cb$
in this order. Note that if $c$ is in $RR$ and $d$ is in $RB$ 
then
$ca$ and $db$ are coloured red and blue, whereas
$cb$ and $da$ are are both red, so $abcd$ is improperly coloured.
We deduce that one of $RR$ and $RB$ has size at most $\epsilon n/3$, since otherwise
we would have at least $\epsilon^2 n^2/9$ improperly colored $4$-tuples containing $ab$.
The same argument applies when take one point from each of $BB$ and $RB$, or
$RR$ and $BR$, or $BB$ and $BR$. Therefore, either $RB$ and $BR$ each have
size at most $\epsilon n/3$, or $RR$ and $BB$ each have size at most $\epsilon n/3$.

In the case when $RB$ and $BR$ each have size at most $\epsilon n/3$ we look
at the pairs in $RR \cup BB$. If $c$ and $d$ are both in $RR$ then both of the
opposite pairs $\{ac,bd\}$ and $\{ad,bc\}$ are coloured red. If $cd$ is
coloured blue then $abcd$ is improperly coloured, so all
but at most $\epsilon^2 n^2/10$ pairs in $RR$ are coloured red. Similarly all but
at most $\epsilon^2 n^2/10$ pairs in $BB$ are coloured red, and all but
at most $\epsilon^2 n^2/10$ pairs with one vertex in $RR$ and one in $BB$ are
coloured blue. Define a partition $V = V_1 \cup V_2$, where $V_1$ contains
$RR$, $V_2$ contains $BB$ and the remaining vertices are distributed
arbitrarily. 
Note that all the incorrect $4$-tuples 
with respect to this partition belong to the one of the following three groups.
 
(i)\, Improperly colored $4$-tuples. There are at most $\epsilon^2 
n^4/130$ 
of those.

(ii)\, Properly colored $4$-tuples which use at least one vertex in $RB 
\cup 
BR$.
There are at most $\big(2\epsilon n/3\big){n \choose 3}$ such $4$-tuples. 

(iii)\, Properly colored $4$-tuples
which contain either a red pair of vertices with one vertex in $RR$ and one 
in $BB$, or contain a blue pair of vertices from $RR$ or from $BB$. There at most
$\big(3\epsilon^2 n^2/10\big) {n \choose 2}$ such $4$-tuples.
 
Therefore all but at most 
$\frac{\epsilon^2 n^4}{130}+2\frac{\epsilon n}{3}{n \choose 3}+3\frac{\epsilon^2 n^2}{10} {n 
\choose 2}< \epsilon n^4$
$4$-tuples are correct with respect to this partition. 

The case when  $RR$ and $BB$ each have size at most $\epsilon n/3$ can be treated 
similarly. Here the conclusion is that all but at most $\epsilon^2 n^2/5$ pairs 
within $RB$ or $BR$ are coloured blue, and all but at most  $\epsilon^2 n^2/10$ pairs 
with one
vertex in $RB$ and one in $BR$ are coloured red. Then, similarly as above one can show that with
respect to a partition where
$V_1$ contains $RB$, $V_2$ contains $BR$ and the remaining vertices are distributed
arbitrarily, all but at most $\epsilon n^4$ $4$-tuples are correct.
\hfill $\Box$

Using the stability theorem we can now prove the following exact Tur\'an result.

\begin{theo}
Let $H$ be a $4$-uniform hypergraph on $n$ vertices that does not contain a 
copy of ${\cal C}^{(4)}_3$ and let $n$ be sufficiently large. Then the number of 
edges in $H$ is at most $b(n)$, with equality only when $H$ is one of
at most $2$ hypergraphs ${\cal B}_n$.
\end{theo}

\noindent {\bf Proof.}\
Let $H$ be a $4$-uniform hypergraph on $n$ vertices, which has
$e(H) \geq b(n)$ and contains no ${\cal C}^{(4)}_3$. 
First 
we claim that we can assume that $H$ has minimum degree at least
$b(n)-b(n-1)$. Indeed, suppose that
we have proved the result under this assumption for all $n \geq n_0$.
Construct a sequence of hypergraphs $H=H_n, H_{n-1},\cdots$ where 
$H_{m-1}$ is obtained from
$H_{m}$ by deleting a vertex of degree less than $b(m)-b(m-1)$. 
By setting $f(m)=e(H_m)-b(m)$ we have
$f(n) \geq 0$ and $f(m) \geq f(m+1)+1$.
If we can continue this process to obtain a hypergraph $H_{n_0}$ then 
$n-n_0 \leq \sum_{m=n_0}^{n-1}\big(f(m)-f(m+1)\big)\leq f(n_0) \leq 
{n_0 \choose 4}$, which is a contradiction for $n$ sufficiently large. Otherwise we 
obtain a hypergraph $H_{n'}$ with $n>n'>n_0$ having minimal 
degree at least $b(n')-b(n'-1)$ 
and without a ${\cal C}^{(4)}_3$. Then by the above assumption 
$e(H_{n'})\leq 
b(n')$ and again we obtain a contradiction, since
$$e(H)=e(H_n) \leq b(n')+\sum_{n'<m\leq n} \big(b(m)-b(m-1)-1 
\big)< b(n)\,.$$
Substituting from equation (\ref{diff}) 
we can assume $H$ has minimum degree
\begin{equation}
\label{eqn5}
\delta(H) \geq b(n) - b(n-1) > \frac{1}{12}n^3 -\frac{1}{2} n^2\,.
\end{equation}
 
Given a partition of $V(H) = V_1 \cup V_2$, we call an edge $abcd$ of $H$
{\em good} if $abcd$ is a good $4$-tuple (as defined before) with respect 
to this partition; otherwise we call it {\em bad}.
By Theorem \ref{stability} there is a partition 
with all but at most
$10^{-25} n^4$ edges of $H$ being good. Let $V(H) = V_1 \cup V_2$
be the partition which 
minimises the number of bad edges. With respect to this partition, every 
vertex belongs to at least as many good edges as bad edges, or we can 
move it to the other class of the partition.
Also, by definition, there are at most 
$b(n)$ good $4$-tuples with respect to any partition. 
We must have $\big| |V_1| - n/2 \big| < 10^{-6} n$ and
$\big| |V_2| - n/2 \big| < 10^{-6} n$. Otherwise by equation
(\ref{heq}) we get
$$e(H) < \frac{1}{48}\Big( \big(n^2-3n+4\big)^2 -
    \big(4 \cdot 10^{-12} n^2-3n+4 \big)^2 \Big) + 10^{-25} n^4 < b(n),$$
which is a contradiction.

Note that there is no
pair of vertices $ab$ for which there are both $10^{-10} n^2$ pairs $cd$ such
that $abcd$ is a good edge and $10^{-10} n^2$ pairs $ef$ such
that $abef$ is an bad edge. Indeed, each such $cd$ and  
$ef$ which are disjoint 
give a $4$-tuple $cdef$ which is good, but cannot be an edge as it would
create a ${\cal C}^{(4)}_3$. Moreover, every $4$-tuple can be obtained 
at most $3$ times in this way, and every $cd$ is disjoint from all but at 
most $2n$ pairs $ef$. Thus at least $10^{-10} n^2\big(10^{-10} 
n^2-2n\big)/3>10^{-21} n^4$ good $4$-tuples  are not edges of $H$, and 
therefore $e(H) < b(n) - 10^{-21} n^4 + 10^{-25} n^4 < b(n)$,
which is a contradiction.

The next step of the proof is the following claim.

\begin{claim} 
\label{cl1}
Any vertex of $H$ is contained in at most $10^{-5} n^3$
bad edges.
\end{claim}

\noindent {\bf Proof.}\ Suppose some vertex $a$ belongs to 
$10^{-5} n^3$
bad edges. Call another vertex $b$ {\em good} if there
are at most $10^{-10} n^2$ pairs $cd$ such that $abcd$ is a bad edge,
otherwise call $b$ {\em bad}. By the above discussion, for every bad
vertex $b$ there are at most $10^{-10} n^2$ pairs $ef$ such
that $abef$ is a good edge. 
Note that there are at least $10^{-5} n$
bad vertices, otherwise we would only have at most
$10^{-5} n \cdot {n \choose 2} + (1-10^{-5})n \cdot 10^{-10} n^2 < 10^{-5} n^3$
bad edges through $a$, which is contrary to our assumption.
By choice of partition there are
at least as many good edges containing $a$ as bad.
We know that $a$ has degree at least
$\frac{1}{12}n^3 -\frac{1}{2} n^2$, at least half of which is good,
so there are at least $n/24$ good vertices.
 
Suppose that the number of good vertices is $\alpha n$, and
so there are $(1-\alpha)n-1$ bad vertices.
We can count the edges containing $a$ as follows.
By definition there are at most $10^{-10} n^3$ such good edges containing 
a bad vertex,
and at most $10^{-10} n^3$ such bad edges containing a good vertex.
Now we bound the number of remaining good edges. Note that these edges 
only contain
good vertices. Looking at the vertices of such an edge in some order, we 
can select the
first $2$ vertices in $\alpha n(\alpha n -1)$ ways. Since the edge is 
good, the choice of $2$ vertices together with $a$ restricts the fourth 
vertex to lie in some particular 
class $V_i$, so it 
can be chosen in at most $\big(\frac{1}{2}+10^{-6}\big)n$ ways. Note 
that we have 
counted each
edge $6$ times, so we get at most 
$\alpha n(\alpha n -1)\big(\frac{1}{2}+10^{-6}\big)n/6<
\big(\alpha^2 + \frac{1}{2} \cdot 10^{-5}\big) \frac{1}{2}{n \choose 3}$
edges. Similarly there are at most 
$\big((1-\alpha)^2 + \frac{1}{2} \cdot 10^{-5} \big) \frac{1}{2}{n \choose 
3}$ remaining bad edges through $a$. Since $1/24 \leq \alpha \leq 
1-10^{-5}$, 
in total the number of edges containing $a$ is bounded by 
$$\left(\alpha^2 + \frac{1}{2} \cdot 10^{-5} \right) \frac{1}{2}{n \choose 
3}+
\left((1-\alpha)^2 + \frac{1}{2} \cdot 10^{-5} \right) \frac{1}{2}{n 
\choose 3}+
2\cdot 10^{-10}n^3<\frac{1}{12}n^3 -\frac{1}{2} n^2 < \delta(H)\,.$$
This contradiction proves the claim.
\hfill $\Box$ 

\vspace{0.15cm}
Now write $|V_1|=n/2+t$, $|V_2|=n/2-t$ with $-10^{-6} n < t < 10^{-6} n$. By
possibly renaming the classes (i.e. replacing $t$ with $-t$)
we can assume that $d(n,t)<d(n,-t)$. Then any vertex of $V_1$ belongs to
$d(n,t)$ good $4$-tuples. Now $d(n,t)$ is the minimum degree
of ${\cal B}(n,t)$, which is certainly at most the maximum degree of 
${\cal B}_n$.
Comparing with equation (\ref{dmax}) we see that any vertex
of $V_1$ belongs to at most
$\frac{1}{12}n^3 - \frac{1}{2} n^2 + n^{3/2}$ good $4$-tuples. From
now on this will be the only property of $V_1$ we use that might possibly not
be a property of $V_2$. We will eventually end up showing the
same bound on the number of good $4$-tuples containing a vertex of $V_2$.
Then the whole argument will apply verbatim switching $V_1$ for $V_2$.

We will use this property in the following manner. Suppose $a$ is
a vertex of $V_1$ for which $K$ of the good $4$-tuples containing $a$ are not
edges of $H$. Then there are at most
$\frac{1}{12}n^3 - \frac{1}{2} n^2 + n^{3/2} - K$ good edges containing $a$,
so by (\ref{eqn5}) there must be at least 
$$\delta(H)-\left(\frac{1}{12}n^3 - \frac{1}{2} n^2 + n^{3/2} - 
K\right)\geq\left(\frac{1}{12}n^3 - \frac{1}{2} n^2\right)-
\left(\frac{1}{12}n^3 - \frac{1}{2} n^2 + n^{3/2} -K\right)
= K-n^{3/2}$$ bad edges containing $a$.
Similarly, if $a'$ is a vertex in $V_2$ then it belongs to at most
$$d(n,-t) = \frac{1}{12}(n^3-6n^2+8n+12t^2) + \frac{1}{12}|6tn-8t^3-16t|
< \frac{1}{12}n^3 - \frac{1}{2} n^2 + 10^{-6} n^3$$
good edges. Thus, if it belongs to $L$
good $4$-tuples which are not edges of $H$ then it must belong to at least 
$L-10^{-6}n^3$ bad edges.

Suppose for the sake of contradiction that there is some bad edge incident
with $V_1$. Denote the set of bad edges containing some vertex 
$v$ by ${\cal Z}(v)$.
Let $a$ be a vertex in $V_1$ belonging to the maximum number of bad edges
and let $Z=|{\cal Z}(a)|$. Note that $Z>0$.
For every bad edge $abcd$ containing $a$, consider
a partition of its vertices into pairs, say $ac$ and $bd$. 
Recall that there are $2$ types of pairs, one
type consisting of pairs with both vertices in one class, the other 
consisting of pairs that have one point of each class.
By definition of a bad edge, $ac$ and $bd$ are pairs
of the same type. If $ef$ is any pair of the other
type which is disjoint from both of them,
then $acef$ and $bdef$ are good $4$-tuples. One of them is not
an edge of $H$, or we get a ${\cal C}^{(4)}_3$. The number of such pairs 
$ef$ is clearly at least
$$\min\bigg\{\Big(|V_1|-4\Big)\Big(|V_2|-4\Big), {|V_1|-4 \choose 
2}+{|V_2|-4 \choose 2}\bigg\} 
\geq \left( \frac{1}{4} - 10^{-12} \right) n^2 - O(n)> n^2/5\,.$$
Let ${\cal Z}_1(a)$ be those bad edges for which there is some partition into pairs 
$ac$ and $bd$, so that for at least
$n^2/10$ of the pairs $ef$ defined above, the good $4$-tuple $acef$
is not an edge. Let ${\cal Z}_2(a)={\cal Z}(a)-{\cal Z}_1(a)$, and write 
$Z_i = |{\cal Z}_i(a)|$ for $i=1,2$. Then one of $Z_1$,$Z_2$ is at least $Z/2$. 

{\bf Case 1:}
Suppose $Z_1 \geq Z/2$. Let $C$ be the (non-empty) set of vertices $c$
such that there is some edge $abcd$ in ${\cal Z}_1(a)$, and $acef$ 
is a good non-edge for at least $n^2/10$ pairs $ef$.
Then we have at least $|C|n^2/30$ good non-edges containing $a$, as we
count each $acef$ at most $3$ times. This implies that there are at least
$|C|n^2/30 - n^{3/2} \geq |C|n^2/31$ bad edges containing
$a$ and therefore $n^2/31 \leq Z/|C|$. Since every edge in 
${\cal Z}_1(a)$ contains at most $3$ vertices of $C$ there
exists $c\in C$ which is contained in at least
$|{\cal Z}_1(a)|/(3|C|)=Z_1/(3|C|)\geq Z/(6|C|)>n^2/200$ bad edges.
Fix one such $c$. 

Note that a graph with $n$ vertices and $m$ edges contains a matching of 
size at least $m/2n$, since  otherwise there
is a set of fewer than $m/n$ vertices that cover all the edges of the 
graph, which is impossible by direct counting. 
Consider the set of pairs $bd$ such that $abcd$ is a bad edge.
Then there exists a matching $M$ of size at least
$n/400$ so that for each $bd$ in $M$ we have that 
$abcd$ is a bad edge of $H$. Partition such an edge into pairs
$ab$ and $cd$. Then, as we explained above, there are at least $n^2/5$ pairs 
$ef$ such that one of the $4$-tuples $abef$ and $cdef$ is a good 
non-edge. Since $M$ is 
a matching we count each such $4$-tuple at most $3$ times, so
one of $a$ or $c$ belongs to at least 
$\frac{1}{3} \cdot \frac{1}{2} \cdot \frac{n^2}{5} \cdot \frac{n}{400}
= n^3/12000$ good non-edges. Therefore it belongs to
at least $n^3/12000 - 10^{-6} n^3 > 10^{-5} n^3$ bad edges,
which contradicts Claim \ref{cl1}. 

{\bf Case 2:}\
Now suppose $Z_2 \geq Z/2$. Note that every bad edge containing $a$
contains at least one other point of $V_1$, so there is
some $b \in V_1$ belonging to at least $Z_2/n$
edges of ${\cal Z}_2(a)$. Fix one such $b$. Suppose $cd$ is a pair such 
that $abcd$ is in ${\cal Z}_2(a)$, and consider any partition of $abcd$ 
into pairs $p_1$, $p_2$
with $a$ in $p_1$ and $b$ in $p_2$. Then, by definition of ${\cal 
Z}_2(a)$, there
are at least $n^2/10$ pairs $ef$ such that $p_2 \cup ef$ is
a good non-edge. Let $C$ be the set of vertices $c$ for which
there exists a vertex $d$ such that $abcd$ is an edge of ${\cal Z}_2(a)$. 
Then
there are at least $|C|n^2/30$ good non-edges containing $b$, as we 
count each $bcef$ at most $3$ times.
Thus, there are at least $|C|n^2/30 - n^{3/2} > 
|C|n^2/50$ bad edges containing $b$. By maximality of $Z$ we have 
$|C|n^2/50 \leq |{\cal Z}(b)| \leq  Z$.
Note that each edge in ${\cal Z}_2(a)$ that contains $b$ is obtained by
picking a pair of vertices in $C$, so
$Z/(2n)\leq Z_2/n \leq {{|C|} \choose 2} < 1250Z^2/n^4$. 
Therefore $Z \geq n^3/2500$, which again contradicts Claim \ref{cl1}. 

We conclude that there are no bad
edges incident to the vertices of $V_1$, i.e. all bad edges
have all $4$ vertices in $V_2$.
We can use this information to give more precise bounds on
the sizes of $V_1$ and $V_2$. We recall that $|V_1|=n/2+t$,
$|V_2|=n/2-t$ and $d(n,t)<d(n,-t)$.
Suppose that $|t| \geq \sqrt{n}$, so that
$6|t|n-8|t|^3-16|t|<-2n^{3/2}$ and by (\ref{degree})
$$d(n,t) = \frac{1}{12}\big(n^3-6n^2+8n+12t^2\big) + 
\frac{1}{12}\big(6|t|n-8|t|^3-16|t|\big)
< \frac{1}{12}n^3 - \frac{1}{2}n^2 - \frac{n^{3/2}}{12}<\delta(H)\,.$$
This is a contradiction, since the
vertices of $V_1$ only belong to good edges, of which there
are at most $d(n,t)<\delta(H)$. 
Therefore $|t| < \sqrt{n}$. Now we can bound the number of
good $4$-tuples containing a vertex of $V_2$. By 
(\ref{degree}), this number is at most
$$d(n,-t) = \frac{1}{12}\big(n^3-6n^2+8n+12t^2\big) + 
\frac{1}{12}\big(8|t|^3-6|t|n+16|t|\big)
< \frac{1}{12}n^3 - \frac{1}{2}n^2 + n^{3/2}\,.$$
Now the same argument as we used to show that no bad edges are
incident with the vertices of $V_1$ shows that none are incident with 
$V_2$ either. We
conclude that all edges are good. Then by definition of $b(n)$
we have $e(H) \leq b(n)$, with equality only when $H$ is a ${\cal B}_n$,
so the theorem is proved. \hfill $\Box$

\section{Proof of Frankl's conjecture}

In this section we will prove the general case of the Frankl conjecture.
We recall that ${\cal C}^{(2k)}_3$ is the $2k$-uniform hypergraph 
with three edges $\{P_1\cup P_2,P_2\cup P_3,P_3\cup P_1\}$, where 
$P_1,P_2,P_3$
are pairwise disjoint sets of $k$ vertices. We can obtain a large
${\cal C}^{(2k)}_3$-free graph on $n$ vertices by partitioning an 
$n$-element set $V$ into $2$ parts $V_1,V_2$ and taking those edges which
intersect each part $V_i$ in an odd number of elements. To see this, consider
any $P_1,P_2,P_3$ that are pairwise disjoint sets of $k$ vertices.
Then $|V_1 \cap P_i|$ and $|V_1 \cap P_j|$ have the same parity for
some pair $ij$, so $P_i\cup P_j$ is not an edge.

Note that this construction is the same as the one we described for ${\cal C}^{(4)}_3$
when $k=2$. In the $4$-uniform case we were able to calculate the sizes
of the parts that maximise the number of edges. For general $k$ this
is an interesting problem in enumerative combinatorics, that is equivalent
to finding the minima of binary Krawtchouk polynomials. 
These polynomials play an important r\^ole in the analysis
of binary Hamming association schemes and so many of
their properties are well-known in this context (see, e.g., \cite{KL}). In 
particular, the location
of their roots is an important problem, but we will need here only a crude 
estimate that follows easily from known results. In the first subsection of 
this section we
will state this estimate and apply it to various parameters of our
construction. The rest of the proof follows the same
broad outline as that of the $4$-uniform case,
in that it falls naturally into two parts.
We will prove the stability part in the second subsection, and
the full result we defer to the final subsection.

\subsection{Binary Krawtchouk polynomials}

Let ${\cal B}^{(2k)}(n,t)$ denote the $2k$-uniform hypergraph obtained 
by partitioning an $n$-element set into two parts
with sizes $\frac{n}{2}+t$ and $\frac{n}{2}-t$, and taking as edges
all $2k$-tuples with odd intersection with each part. Let $b_{2k}(n,t)$ be the number
of edges in ${\cal B}^{(2k)}(n,t)$ and let $d_{2k}(n,t)$ be the degree of any vertex
belonging to the side with size $\frac{n}{2}+t$. Then the
vertices on the side with size $\frac{n}{2}-t$ have degree
$d_{2k}(n,-t)$. 

The {\em binary Krawtchouk polynomials} $K_m^n(x)$
can be defined by the generating function
$$\sum_{m=0}^n K_m^n(x) z^m = (1-z)^x (1+z)^{n-x}.$$
From here we get
the explicit expression $K_m^n(x) = \sum_{i=0}^m (-1)^i {x \choose i}{{n-x} \choose {m-i}}$.
Recall that $b_{2k}(n,t)$ was the number of 
$2k$-tuples with odd intersection with both parts in the above partition of an 
$n$-element set and so ${n \choose {2k}}-b_{2k}(n,t)$ is the number of
$2k$-tuples with even intersection with these parts. 
This implies that 
$\big({n \choose {2k}}-b_{2k}(n,t)\big)-b_{2k}(n,t) = \sum_{i=0}^{2k} (-1)^i 
{{n/2+t} \choose i}{{n/2-t} \choose {2k-i}}
= K_{2k}^n(n/2+t)$, which gives
\begin{equation}
\label{b2k}
b_{2k}(n,t) = \frac{1}{2} \bigg( {n \choose {2k}} - K_{2k}^n(n/2+t) \bigg),
\end{equation}
so maximising $b_{2k}(n,t)$ is equivalent to finding the minimum of $K_{2k}^n(x)$.
Similarly, we can also express the degrees of ${\cal B}^{(2k)}(n,t)$ in terms 
of Krawtchouk polynomials. Indeed, by definition, $d_{2k}(n,t)$ is the number 
of $(2k-1)$-tuples with even intersection with the first part
in the partition of an $(n-1)$-element set in two parts with sizes $n/2+t-1$ 
and $n/2-t$, and therefore ${n-1 \choose {2k-1}}-d_{2k}(n,t)$
is the number of $(2k-1)$-tuples with odd intersection with this part.
Then, 
$d_{2k}(n,t)-\big({n-1 \choose {2k-1}}-d_{2k}(n,t)\big) = \sum_{i=0}^{2k-1} 
(-1)^i {{n/2+t-1} \choose i}{{n/2-t} \choose {2k-1-i}}
= K_{2k-1}^{n-1}(n/2+t-1)$, i.e.
\begin{equation}
\label{d2k}
d_{2k}(n,t) = \frac{1}{2} \bigg( {n-1 \choose {2k-1}} + K_{2k-1}^{n-1}(n/2+t-1) 
\bigg)
\end{equation}

Note that $K_m^n(x)$ is a polynomial of degree $m$. It is known that it has $m$ simple roots,
symmetric with respect to $n/2$. 
The smallest root is given by the following formula obtained by Levenshtein
\cite{Le}:
$$ r = n/2 - \max \bigg( \sum_{i=0}^{m-2} x_i x_{i+1} \sqrt{(i+1)(n-i)} \bigg), 
$$
where the maximum is taken over $x_i$ with $\sum_0^{m-1} x_i^2 = 1$. From
the Cauchy-Schwartz inequality we see that $n/2-r<\sqrt{mn}$. Note that
$K_{2k}^n(0)=K_{2k}^n(n)={n \choose {2k}}>0$, so the minimum of $K_{2k}^n(x)$ occurs in
the range $n/2 \pm \sqrt{2kn}$.

Let $t^*$ be chosen to maximise the number of edges in ${\cal B}^{(2k)}(n,t)$,
and denote any hypergraph obtained in this manner by ${\cal B}^{(2k)}_n$. Note
that $t^*$ may not be unique, but must satisfy $|t^*|<\sqrt{2kn}$. Also, by 
symmetry we can assume that $t^*>0$.
Write $b_{2k}(n)$ for the number of edges in ${\cal B}^{(2k)}_n$.

\begin{lemma}
\label{kraw}
$(i)$\, $K_m^n(n/2+t) = \sum_{i=0}^{m/2} (-1)^{i+m} {{n/2-t} \choose i}{{2t} 
\choose {m-2i}}.$

\noindent
$(ii)$\, If $c>1$ and $0\leq s \leq c\sqrt{n}$ then $\Big| d_{2k}(n,\pm s) 
- \frac{1}{2} {{n-1} \choose 
{2k-1}} \Big| < (10c^2)^k n^{k-1/2}.$

\noindent
$(iii)$\, $\Big| b_{2k}(n) - \frac{1}{2} {n \choose {2k}} \Big| < (20kn)^k$, 
\,\,
$\Big| d_{2k}(n,\pm t^*) - \frac{1}{2} {{n-1} \choose {2k-1}} \Big| < (20k)^k 
n^{k-1/2}.$ 

\noindent
$(iv)$\, If $C>20^k$ then $d_{2k}(n,C\sqrt{n}) < \frac{1}{2} {{n-1} \choose 
{2k-1}} - 20^k n^{k-1/2}$.

\noindent
$(v)$\, $\Big| b_{2k}(n,\epsilon n) - \Big( \frac{1}{2} {n \choose {2k}} - 
\frac{1}{2} {{2\epsilon n} \choose {2k}} \Big)
\Big| < (10\epsilon)^k n^{2k-1}$, \,\,
$\Big| d_{2k}(n,\epsilon n) -  \Big( \frac{1}{2} {{n-1} \choose {2k-1}} - 
\frac{1}{2} {{2\epsilon n-1} \choose {2k-1}} \Big)
\Big| < (10\epsilon)^k n^{2k-2}.$
\end{lemma}

\noindent {\bf Proof.}\,
(i)\, Rewrite the generating function as $\sum_{m=0}^n K_m^n(n/2+t) z^m = 
(1-z^2)^{n/2-t} (1-z)^{2t}$ and expand.

\noindent
(ii)\, Using part (i) with $t=s-1/2$, and applying (\ref{d2k}), 
we get 
\begin{eqnarray*}
\bigg| d_{2k}(n,s) - \frac{1}{2} {{n-1} \choose {2k-1}} \bigg|
&=& \frac{1}{2} \bigg|K_{2k-1}^{n-1}(n/2+s-1) \bigg|\\
&=& \frac{1}{2} \bigg| \sum_{i=0}^{k-1} (-1)^{i+1} {{n/2-s} \choose i} 
{{2s-1} \choose {2k-1-2i}} \bigg|\\
&<& k \cdot(2c\sqrt{n})^{2k-1} < (10c^2)^k n^{k-1/2}.
\end{eqnarray*}
The corresponding inequality for $d_{2k}(n,-s)$ can be obtained similarly.

\noindent
(iii)\, The second statement follows from part (ii) with $c =\sqrt{2k}\geq
t^*/\sqrt{n}$. To prove the first statement, we use (\ref{b2k}), 
part (i) and again the fact that $0<t^* <\sqrt{2kn}$. Altogether they imply
$$\bigg| b_{2k}(n) - \frac{1}{2} {n \choose {2k}} \bigg| = \frac{1}{2} \Big| 
K_{2k}^n(n/2+t^*)\Big |
< \sum_{i=0}^k {{n/2} \choose i}{{2\sqrt{2kn}} \choose {2k-2i}}
< (k+1) \big(2\sqrt{2kn}\big)^{2k} < (20kn)^k.$$

\noindent
(iv)\, By (\ref{d2k}) we have
\begin{eqnarray*}
d_{2k}(n,C\sqrt{n}) - \frac{1}{2} {{n-1} \choose {2k-1}}
&=& \frac{1}{2} K_{2k-1}^{n-1}\big(n/2+C\sqrt{n}-1\big)
=  \frac{1}{2} \sum_{i=0}^{k-1} (-1)^{i+1} {{\frac{n}{2}-C\sqrt{n}} \choose i} 
{{2C\sqrt{n}-1} \choose {2k-1-2i}}\\
&<& - \frac{1}{2} {{2C\sqrt{n}-1} \choose {2k-1}} +\frac{1}{2} (k-1) (n/2) 
{{2C\sqrt{n}-1} \choose {2k-3}}\\
&=& \big(1+o(1)\big) \left( - \frac{(2C\sqrt{n})^{2k-1} }{2(2k-1)!} 
       + \frac{(k-1)n}{4} \frac{(2C\sqrt{n})^{2k-3}}{(2k-3)!} \right)\\
&<&  -\left( \frac{C^2}{(2k-1)^2} - \frac{k-1}{2} 
\right) \frac{(2C)^{2k-3}}{(2k-3)!} n^{k-1/2}
< -  20^k n^{k-1/2}.
\end{eqnarray*}

\noindent
(v)\, Using formula for $K_{2k}^n(n/2+t)$ from part (i) together with (\ref{b2k}) we obtain that
\begin{eqnarray*}
\bigg| b_{2k}(n,\epsilon n) - \Big( \frac{1}{2} {n \choose {2k}} - \frac{1}{2} {{2\epsilon n} 
\choose {2k}} \Big)\bigg| &=& \bigg| \frac{1}{2} K_{2k}^n(n/2+\epsilon n) - \frac{1}{2} 
{{2\epsilon n} \choose {2k}} \bigg|\\&=&
\frac{1}{2} \bigg| \sum_{i=1}^{k} (-1)^i {{n/2-\epsilon n} \choose i}{{2\epsilon n} \choose 
{2k-2i}} \bigg|\\
&<& \frac{n}{2} \cdot (2\epsilon n)^{2k-2}+O\big(n^{2k-2}\big) < (10\epsilon)^k n^{2k-1}.
\end{eqnarray*}
The proof of the inequality for $d_{2k}(n,\epsilon n)$ can be obtained similarly and we omit it 
here. \hfill $\Box$

We remark that these simple estimates are sufficient for our purposes, but the
location of the roots and asymptotic values for Krawtchouk polynomials in the oscillatory
region are known with more precision (see, e.g., \cite{KL}). With this 
information one could find better
estimates for $b_{2k}(n)$, and possibly how many different choices of $t$ give the
maximum number of edges.

We conclude this section with an estimate on the difference of successive values
of $b_{2k}(n)$.

\begin{lemma}
\label{difference}
$b_{2k}(n)-b_{2k}(n-1) \geq \frac{1}{2} {n-1 \choose {2k-1}}$
\end{lemma}

\noindent {\bf Proof.}\,
Suppose that ${\cal H} = {\cal B}^{(2k)}(n-1)$ has $b_{2k}(n-1)$ edges and has 
parts $V({\cal H})=A\cup B$.
Let ${\cal H}_1$ be obtained from ${\cal H}$ by adding
a vertex $v_1$ to $A$, together with all the $2k$-tuples containing $v_1$ and 
having odd intersections with $A\cup v_1$ and $B$.
Let ${\cal H}_2$ be similarly obtained by adding a vertex $v_2$ to $B$,
together with corresponding edges.
By definition each ${\cal H}_i$ has at most $b_{2k}(n)$ edges, so the degree  
of each $v_i$ is
a lower bound for $b_{2k}(n)-b_{2k}(n-1)$. On the other hand, for each 
$(2k-1)$-tuple $X$ of vertices in ${\cal H}$ there is exactly one $i$ such that 
$X \cup v_i$ is
an edge of ${\cal H}_i$, so one of the $v_i$ has degree
at least $\frac{1}{2} {n-1 \choose {2k-1}}$. \hfill $\Box$ 

\subsection{A stability result for ${\cal C}^{(2k)}_3$}

In this subsection we prove a stability result for ${\cal C}^{(2k)}_3$.
We start by recalling a version of the Kruskal-Katona
theorem due to Lovasz. Write $[m]=\{1,\cdots,m\}$, let 
$[m]^{(k)}$
denote the subsets of $[m]$ of size $k$, and suppose ${\cal A} \subset [m]^{(k)}$.
The {\em shadow} of ${\cal A}$ is $\partial {\cal A} \subset [m]^{(k-1)} $ consisting
of all sets of size $k-1$ that are contained in some element of ${\cal A}$.
For any real $x$ write ${x \choose k} = x(x-1)\cdots(x-k+1)/k!$. 
The following result appears in \cite{Lo} (Exercise 13.31). 

\begin{prop}
\label{kk}
If ${\cal A} \subset [m]^{(k)}$ and $|{\cal A}| = {x \choose k}$
then $|\partial {\cal A}| \geq {x \choose {k-1}}$. \hfill $\Box$
\end{prop} 
     
Suppose we have a $2k$-uniform hypergraph $H$ 
and a partition of the vertex set $V(H)=V_1 \cup V_2$.
Our terminology for $2k$-tuples matches that of the $4$-uniform case.  
We call a $2k$-tuple of vertices {\em good} if it
intersects each $V_i$ in an odd number of elements; otherwise we call it 
{\em bad}. We call a $2k$-tuple {\em correct} if it is either a good edge 
or a bad non-edge; otherwise we call it {\em incorrect}.

\begin{theo} 
\label{fullstab}
For every $\epsilon>0$ there is $\eta>0$ so that if $H$ is a 
${\cal C}^{(2k)}_3$-free
$2k$-uniform hypergraph with $e(H)>\frac{1}{2}{n \choose {2k}} - \eta n^{2k}$ then there
is a partition of the vertex set as $V(H)=V_1 \cup V_2$ such that all but
$\epsilon n^{2k}$ $2k$-tuples are correct.
\end{theo}

\noindent {\bf Proof.}\,
Define an auxiliary graph $G$ whose vertices are all $k$-tuples of vertices of 
$H$, and where the $k$-tuples $P_1$ and  $P_2$ are adjacent exactly when 
$P_1\cup P_2$ is 
an edge of $H$. Since $H$ is ${\cal C}^{(2k)}_3$-free we see that $G$ is 
triangle-free.
Also, each edge of $H$ creates exactly $\frac{1}{2}{{2k} \choose k}$ edges in $G$ (corresponding
to the ways of breaking a $2k$-tuple into two $k$-tuples) so
$$e(G) > \frac{1}{2}{{2k} \choose k} \left( \frac{1}{2}{n \choose {2k}} - 
\eta n^{2k} \right)
       > \Big(1-(k!)^22^{2k}\eta\Big) \frac{1}{2} {{n \choose k} \choose 2}.$$

Choose $\eta$ so that the Simonovits stability theorem (see Section 2) applies 
with $\eta'=(k!)^22^{2k}\eta$,
$N={n \choose k}$ and $\epsilon'=10^{-6k^2}\epsilon^k$. We can
also require that $\eta < 10^{-6k^2}\epsilon^k$. We get a partition of the
$k$-tuples of vertices of $H$ as $U_0 \cup U_1$, where all but
$\epsilon' N^2=\epsilon'{n \choose k}^2 < 10^{-6k^2}\epsilon^k n^{2k}$ edges of $H$ are 
formed by taking
a $k$-tuple from $U_0$ and a $k$-tuple from $U_1$. 

We will think of the sets $U_i$ as determining a $2$-coloring of all
$k$-tuples, and say that the $k$-tuples in $U_i$ have colour $i$.
A $2k$-tuple $I$ will be called {\em properly coloured} if, either it is 
an edge of $H$ and 
however we partition $I$ into $k$-tuples $P_1$ and $P_2$ they have 
different colours, or it is not an edge of $H$ and for any partition of 
$I$ into two $k$-tuples they have the same color.

An improperly coloured $2k$-tuple is either an edge that is the union of two
$k$-tuples of the same colour or
a non-edge which is the union of two $k$-tuples with different colours.
There are at most $10^{-6k^2}\epsilon^k n^{2k}$ of the former 
$2k$-tuples,
and the number of latter is at most
$$|U_1||U_2| - \Big(e(G)-\epsilon'N^2 \Big)\leq 
\frac{(k!)^22^{2k}\eta}{2}\,\frac{N^2}{2} +\epsilon'N^2 \leq 
\left(\frac{(k!)^22^{2k}\eta}{4(k!)^2}+\frac{\epsilon'}{(k!)^2}\right)n^{2k}
\leq 10^{-5k^2-1}\epsilon^k
n^{2k}\,.$$
 Therefore all but 
$\big(10^{-6k^2}\epsilon^k +10^{-5k^2-1}\big)\epsilon^k n^{2k}<10^{-5k^2}\epsilon^k 
n^{2k}$ $2k$-tuples 
are properly coloured.

A simple counting argument shows that there is a $k$-tuple $P$
so that for all but ${{2k} \choose k}10^{-5k^2}\epsilon^k n^{2k}/{n \choose k}$ $< 
10^{-4k^2}\epsilon^k n^k$
other $k$-tuples $Q$ the $2k$-tuple $P\cup Q$ is properly coloured. 
Without loss of generality $P$ has colour $0$. We will call
a $k$-tuple $Q$ {\em proper} if $P\cup Q$ is properly coloured;
otherwise it is {\em improper}. Then by definition there
are at most $10^{-4k^2}\epsilon^k n^k$ improper $k$-tuples.
Call a $(k-1)$-tuple $X \subset V-P$ {\em abnormal} if there are at
least $2^{-3k}\epsilon n$ vertices $x \in V-(P \cup X)$ for which $X \cup 
x$ is
improper; otherwise call it {\em normal}.
It is easy to see that there are at most 
$k\cdot 10^{-4k^2}\epsilon^k n^k/(2^{-3k}\epsilon n)<10^{-3k^2}\epsilon^{k-1} 
n^{k-1}$ abnormal $(k-1)$-tuples.

We partition the vertices of $V- P$ 
according to the colour of the $k$-tuples that they form when they
replace an element of $P$. To be precise, we fix an order
$p_1,\cdots,p_k$ of $P$ and partition into $2^k$ parts $V-P=\bigcup V_{\bf 
s}$, where
${\bf s}=(s_1, \ldots, s_k) \in \{0,1\}^k$ and a vertex $x$ belongs to 
$V_{\bf s}$
iff $(P-p_i)\cup x$ has colour $s_i$ for every $1 \leq i 
\leq k$.

Consider a $(k-1)$-tuple $X=x_1\cdots x_{k-1}$ and suppose $a$ is a 
vertex such that $X\cup a$ is
proper. Fix $1 \leq i \leq k$ and consider the partitions
$P\cup X\cup a = (P)\bigcup (X \cup a) = \big((P-p_i)\cup a\big)\bigcup 
(X\cup p_i)$. 
Let $V_{\bf s}$ be the class
containing $a$, so that $(P-p_i)\cup a$ has colour $s_i$. We recall
that $P$ has colour $0$, so if also $s_i=0$ then to be
properly coloured $X\cup a$ must have the same colour as $X \cup p_i$.
On the other hand, if $s_i=1$ then $X\cup a$ and $X\cup p_i$ must have 
different colours.
If we write $c_X(v)$ for the colour of $X\cup v$ for any vertex $v$, then this
can be summarised as 
\begin{equation}
\label{eqn1}
\text{If} \quad a \in V_s \quad \text{and} \quad X\cup a \quad
\text{is proper, then} \quad c_X(a)+s_i=c_X(p_i)~~ 
(\hspace{-0.4cm}\mod 2)\,.
\end{equation}

Suppose there are $2$ classes $V_{\bf s}$ and $V_{\bf s'}$ both of
size at least $2^{-2k}\epsilon n$. Since ${{2^{-2k}\epsilon n} \choose 
k-1} > 10^{-3k^2}\epsilon^{k-1} n^{k-1}$
some $(k-1)$-tuple $X \subset V_s$ is normal. This means that there
are at most $2^{-3k}\epsilon n$ vertices $x \in V-(P\cup X)$ for which 
$X\cup x$ is
improper, so there is $a \in V_{\bf s}$ and $b \in V_{\bf s'}$
such that $X\cup a$ and $X\cup b$ are proper. 
For any pair of indices $i,j$ we have $c_X(a)+s_i=c_X(p_i)$, 
$c_X(a)+s_j=c_X(p_j)$,
$c_X(b)+s'_i=c_X(p_i)$ and $c_X(b)+s'_j=c_X(p_j)$. Adding these
equations gives $s_i+s_j+s'_i+s'_j=0$. If ${\bf s}$ and ${\bf s'}$
differ in some co-ordinate $i$ then this equation shows
that they must also differ in any other co-ordinate $j$.
In other words, if ${\bf s'} \neq {\bf s}$ we must have ${\bf s'} = 
\overline{{\bf s}}$,
where $\overline{{\bf s}}$ denotes the sequence whose $i$th entry
is $1-s_i$. 

Let $V_{\bf s}$ be the largest class, and write $m=|V_{\bf s}|$. Clearly 
$m \geq 2^{-k}(n-k)$.
Then all other classes, except possibly $V_{\overline{{\bf s}}}$, have
size at most $2^{-2k}\epsilon n$. Let ${\cal A}_i$ be the set of proper 
$k$-tuples
contained in $V_{\bf s}$ that have colour $i$. Then
$|{\cal A}_0|+|{\cal A}_1| > {m \choose k} - 10^{-4k^2}\epsilon^k n^k > 
(1-10^{-3k^2}\epsilon) {m 
\choose k}$. Write
$|{\cal A}_i| = \alpha_i {m \choose k}$, so that $\alpha_0+ 
\alpha_1>1-10^{-3k^2}\epsilon$. Suppose both $\alpha_i$ are at least 
$10^{-2k^2}\epsilon$.
Observe that $|{\cal A}_i| = {{\alpha_i^{1/k}m} \choose k} + O(m^{k-1})$, so
by Proposition \ref{kk} we have 
$$|\partial {\cal A}_i| \geq {{\alpha_i^{1/k} m} \choose {k-1}} + 
O(m^{k-2})= \alpha_i^{(k-1)/k} {m \choose {k-1}} + O(m^{k-2}).$$
Note that if $z \leq 2^{-k}$ we have that $z^{-1/k}\geq 2$ and therefore
\begin{eqnarray*}
z^{(k-1)/k} + (1-10^{-3k^2}\epsilon-z)^{(k-1)/k}
&\geq& z^{(k-1)/k} + \big(1-10^{-3k^2}\epsilon-z\big)\\
&\geq& 2z + 1-10^{-3k^2}\epsilon - z = 1+z-10^{-3k^2}\epsilon .
\end{eqnarray*}
Since $z^{(k-1)/k} + (1-10^{-3k^2}\epsilon-z)^{(k-1)/k}$ is concave we have
$\alpha_0^{(k-1)/k} + \alpha_1^{(k-1)/k} \geq 1 +  
10^{-2k^2}\epsilon - 10^{-3k^2}\epsilon \geq 1+10^{-3k^2}\epsilon$.
We deduce that $|\partial {\cal A}_0 \cap \partial {\cal A}_1| > 0$, i.e.
there is a $(k-1)$-tuple $X$ and points $a_0, a_1$ such that
$X\cup a_i$ is proper, with $c_X(a_i)=i$. But equation (\ref{eqn1})
gives $i+s_1=c_X(a_i)+s_1=c_X(p_1)$, for $i=0,1$, which is
a contradiction. We conclude that there is $t \in \{0,1\}$ for
which $\alpha_{1-t} < 10^{-2k^2}\epsilon$, and so all
but at most $10^{-2k^2}\epsilon {m \choose k} + 
10^{-4k^2}\epsilon^k n^k < 
10^{-2k^2}\epsilon n^k$ $k$-tuples
inside $V_{\bf s}$ have the same colour $t$.

For $0 \leq i \leq k$ let ${\cal D}_i$ be all $k$-tuples with $i$
points in $V_{\overline{\bf s}}$ and $k-i$ points in $V_{\bf s}$ and  
let $\theta_i = 10^{-2k^2}\big( 2k2^{2k} \big)^i\epsilon$.
We claim that for each $i$ all but at most 
$\theta_i n^k$ $k$-tuples 
of ${\cal D}_i$ have colour $t+i$ (mod $2)$. Otherwise, choose
the smallest $i$ for which this is not true. 
By the above discussion $i>0$, and there are at least $\theta_i n^k$ $k$-tuples 
in ${\cal D}_i$ with color $1-(t+i)=t+i-1$ (mod $2)$. Since $i$ was the smallest 
such index all but at most 
$\theta_{i-1} n^k$ $k$-tuples of ${\cal D}_{i-1}$ have colour $t+i-1$ (mod 
$2$). Let $E_{i-1}$ be
the $(k-1)$-tuples $Y$ with $i-1$ points in $V_{\overline{\bf s}}$ and 
$k-i$ points
in $V_{\bf s}$ for which there are at least $2^{-2k} n$ points $y \in 
V_{\bf s}$
such that $Y\cup y$ does not have colour $t+i-1$ (mod $2$).
Then $|E_{i-1}| \leq k \theta_{i-1} n^{k}/(2^{-2k} n)=\frac{1}{2} \theta_i
n^{k-1}$, 
so at most $\frac{1}{2} \theta_i n^k$ $k$-tuples contain an element
of $E_{i-1}$. Recall that there are at most $10^{-4k^2}\epsilon^kn^k$ 
improper $k$-tuples
and at most $10^{-3k^2}\epsilon^{k-1} n^{k-1} \cdot n=
10^{-3k^2}\epsilon^{k-1} n^k$ $k$-tuples that contain some abnormal 
$(k-1)$-tuple.
Since $10^{-4k^2}\epsilon^k+10^{-3k^2}\epsilon^{k-1} < 10^{-2k^2-1}\epsilon<\theta_i/2$ we 
can find a proper $k$-tuple 
$K\in {\cal D}_i$
such that $K$ has color $t+i-1$ (mod $2$) and
for any $(k-1)$-tuple $Y \subset K$ we have $Y$ normal
and $Y \notin E_{i-1}$.

Since $i>0$, there is  $x \in K \cap V_{\overline{\bf s}}$. Let $Y = K-x$. Since $Y$ is normal
there are at most $2^{-3k}\epsilon n$ vertices $y$ such that $Y\cup y$ is improper,
and by definition of $E_{i-1}$ there are at most $2^{-2k} n$ points 
$y \in V_{\bf s}$ such that $Y\cup y$ does not have colour $t+i-1$ (mod $2$). 
Since
$2^{-3k}\epsilon n + 2^{-2k}n < 2^{-k}(n-k)$ there is $y \in V_{\bf s}$
such that $Y\cup y$ is proper and has colour $t+i-1$ (mod $2$). 
Then $c_Y(x)=c_Y(y)=t+i-1$ (mod $2$). But $x \in V_{\overline{\bf s}}$ and 
$y \in V_{\bf s}$, 
so $c_Y(x)+1-s_1=c_Y(p_1)$ and $c_Y(y)+s_1=c_Y(p_1)$, both mod $2$. 
This is a contradiction, so we conclude that all
but at most $\theta_i n^k$ $k$-tuples
of ${\cal D}_i$ have colour $t+i$.

Now partition $V$ into $2$ classes $V_1$, $V_2$ so that
$V_{\bf s} \subset V_1$, $V_{\overline{\bf s}} \subset V_2$, and
the other vertices are distributed arbitarily. Incorrect $2k$-tuples
with respect to this partition belong to the one of the following three groups.

(i)\, Improperly colored $2k$-tuples. There are at most
$10^{-5k^2}\epsilon^kn^{2k}$ of those.

(ii)\, Properly colored $2k$-tuples which use at least one vertex not in
$V_{\bf s} \cup V_{\overline{\bf s}}$. There are at most
$2^k 2^{-2k}\epsilon n {n \choose {2k-1}} < 2^{-k}\epsilon n^{2k}$ such $2k$-tuples.

(iii)\, Properly colored $2k$-tuples which contain a $k$-tuple of ${\cal D}_i$ 
with colour $t+i-1$ (mod $2$). There are at most
$\sum_{i=0}^k \theta_i n^{k}{n \choose k}<\theta_k n^{2k}=
10^{-2k^2}\big(2k2^{2k}\big)^k\epsilon n^{2k}<10^{-k^2}\epsilon n^{2k}$ such $2k$-tuples.

Therefore all but at most $\big(10^{-5k^2}\epsilon^k+2^{-k}\epsilon+10^{-k^2}\epsilon\big)n^{2k}<
\epsilon n^{2k}$ $2k$-tuples are correct
with respect to this partition. This completes the proof of the theorem. \hfill $\Box$

\subsection{The Tur\'an number of ${\cal C}^{(2k)}_3$}
In this subsection we complete the proof of Frankl's conjecture.

\vspace{0.1cm}
\noindent {\bf Proof of Theorem \ref{frankl}.}\
Let $H$ be a $2k$-uniform hypergraph on $n$ vertices, which has
$e(H) \geq b_{2k}(n)$ and contains no ${\cal C}^{(2k)}_3$. 
By the same argument given in the proof in the case $k=2$ we can
assume that $H$ has minimum degree at least $b_{2k}(n)-b_{2k}(n-1)$. Applying
Lemma \ref{difference} gives
\begin{equation}
\label{mindeg}
\delta(H) \geq \frac{1}{2} {n-1 \choose {2k-1}}
\end{equation}
For convenience of notation we set $\eta = (100k)^{-10^k}$. By Theorem \ref{fullstab} there is a partition 
with all but at most $(\eta/20k)^{2k} n^{2k}$ edges of $H$ being good, i.e.,
they have odd intersection with both parts. Let 
$V(H) = V_1 \cup V_2$
be the partition which  minimises the number of bad edges. Then every 
vertex belongs to at least as many good edges as bad edges, or we can 
move it to the other class of the partition. Recall that, by definition, the 
number of good 
$2k$-tuples with respect to this partition is at most $b_{2k}(n)$.
We must have $\big| |V_1| - n/2 \big| < \frac{1}{10}\eta n$ and
$\big| |V_2| - n/2 \big| < \frac{1}{10}\eta n$. Otherwise by Lemma \ref{kraw},
part (5) 
$$e(H) < \frac{1}{2} {n \choose {2k}} - \frac{1}{2} {{2 \cdot \frac{1}{10}\eta n} \choose {2k}} + (10 \cdot \frac{1}{10}\eta)^k n^{2k-1}
+ (\eta/20k)^{2k} n^{2k} < b_{2k}(n),$$
which is a contradiction.

Note that there is no $k$-tuple of vertices $P$ for which there are both $(10k)^{-k}\eta n^k$ 
$k$-tuples $Q$ such that $P \cup Q$ is a good edge and $(10k)^{-k}\eta n^k$ 
$k$-tuples $R$ such that $P \cup R$ is a bad edge. Indeed, each such $Q$ and  
$R$ which are disjoint give a $2k$-tuple $Q \cup R$ which is good, but cannot 
be an edge as it would
create a ${\cal C}^{(2k)}_3$. Moreover, every $2k$-tuple can be obtained 
at most $\frac{1}{2} {{2k} \choose k}$ times in this way, and every $Q$ is disjoint from all but at 
most $k{n \choose {k-1}}$ $k$-tuples $R$. Thus at least 
$(10k)^{-k}\eta n^k \big( (10k)^{-k}\eta n^k - k{n \choose {k-1}} \big) / \big( \frac{1}{2} {{2k} \choose k} \big)
> 2 (\eta/20k)^{2k} n^{2k}$ good $2k$-tuples are not edges of $H$, and 
therefore $e(H) < b_{2k}(n) - 2 (\eta/20k)^{2k} n^{2k} + (\eta/20k)^{2k} n^{2k} < b_{2k}(n)$,
which is a contradiction.

\begin{claim} 
\label{baddeg}
Any vertex of $H$ is contained in at most $\eta n^{2k-1}$ bad edges.
\end{claim}

\noindent {\bf Proof.}\ Suppose some vertex $a$ belongs to $\eta n^{2k-1}$ bad edges.
Call a $(k-1)$-tuple $X$ {\em good} if there are at most $(10k)^{-k}\eta n^k$ 
$k$-tuples $Q$ such that
$a \cup X \cup Q$ is a bad edge, otherwise call $X$ {\em bad}. By the above 
discussion,
for every bad $(k-1)$-tuple $X$ there are at most $(10k)^{-k}\eta n^k$ $k$-tuples $R$
such that $a \cup X \cup R$ is a good edge. There are at least $\eta n^{k-1}$
bad $(k-1)$-tuples or we would only have 
$\eta n^{k-1} \cdot {n \choose k} + \big({n-1 \choose k-1}-\eta 
n^{k-1}\big) \cdot (10k)^{-k}\eta n^k
< \eta n^{2k-1}$ bad edges through $a$. By choice of partition there are
at least as many good edges containing $a$ as bad. From 
(\ref{mindeg}) we see that $a$ is in at least 
$\frac{1}{4} {{n-1} \choose {2k-1}}$ good edges,
so there are at least $\big(\frac{1}{4} {{n-1} \choose {2k-1}}-
(10k)^{-k}\eta n^{2k-1}\big)/{n-1 \choose k}\geq
n^{k-1}/(2k)!$ good $(k-1)$-tuples.

Suppose there are $\alpha {n \choose {k-1}}$ good $(k-1)$-tuples, where
by the above we see that $(2k)^{-k-1} \leq \alpha \leq 1 - (k-1)! \eta$.
We can count the good edges containing $a$ as follows. By definition 
there are at most ${n \choose {k-1}} \cdot (10k)^{-k}\eta n^k$ such good edges 
containing a bad $(k-1)$-tuple. 
Note that in the remaining good edges every $(k-1)$-tuple is good.
Given any such edge $W$ containing $a$ we
consider ordered triples $(X,Y,b)$, where $X$ and $Y$ are $(k-1)$-tuples, $b$
is a vertex and $X \cup Y \cup b \cup a = W$. Each edge gives
rise to $k{{2k-1} \choose {k-1}}$ such triples. To bound the number of
triples recall that $X$ and $Y$ are good, so can be chosen in at most
$\Big( \alpha {n \choose {k-1}} \Big)^2$ ways. Once
$X$ and $Y$ have been chosen, to make $E$ good $b$ is constrained to lie in some
particular class $V_i$ of the partition, so can be chosen in
at most $\Big( \frac{1}{2} + \frac{1}{10}\eta \Big)n$ ways. This shows that
the number of good edges not containing bad $(k-1)$-tuples is at most
$$ \left( \Big( \alpha {n \choose {k-1}} \Big)^2 
\Big( \frac{1}{2} + \frac{1}{10}\eta \Big)n \right) \Big/ \left( k{{2k-1} 
\choose {k-1}} \right)
< \Big( \alpha^2 + 3 \cdot \frac{1}{10}\eta \Big) \frac{1}{2} {n-1 \choose 
{2k-1}}$$
We can count the bad edges similarly, and deduce that the total number
of edges containing $a$ is at most
$$ \Big( \alpha^2 + (1-\alpha)^2 + 6\cdot\frac{1}{10}\eta \Big)  \frac{1}{2} 
{n-1 \choose {2k-1}}
+ 2 \cdot {n \choose {k-1}} \cdot (10k)^{-k}\eta n^k.$$
From the bounds $(2k)^{-k-1} \leq \alpha \leq 1 - (k-1)! \eta$ we
see that this is at most $\Big( \frac{1}{2} - \eta/2 \Big) {{n-1} \choose 
{2k-1}}$.
This contradicts equation (\ref{mindeg}), so the claim is proved. \hfill $\Box$

\vspace{0.15cm}
Now write $|V_1|=n/2+t$, $|V_2|=n/2-t$ with $-\frac{1}{10}\eta n < t < 
\frac{1}{10}\eta n$. By
possibly renaming the classes (i.e. replacing $t$ with $-t$)
we can assume that $d(n,t)<d(n,-t)$. Now any vertex of $V_1$ belongs to 
$d(n,t)$ good $2k$-tuples,
and $d(n,t)$ is the minimum degree of ${\cal B}^{(2k)}(n,t)$, which is 
certainly at most the maximum degree of 
${\cal B}^{(2k)}_n$. From Lemma \ref{kraw}, part (3) we have a bound 
$d(n,t)<\frac{1}{2} {{n-1} \choose {2k-1}} + (20kn)^{k-1/2}$
but we will only use the weaker bound $d(n,t)<\frac{1}{2} {{n-1} \choose 
{2k-1}} +  10^{4k^2} n^{k-1/2}$.
Later we will show that this weaker bound also holds for $d(n,-t)$, and then the
subsequent argument will apply switching $V_1$ and $V_2$.

\begin{claim}
\label{goodtobad}
1. If $a$ is a vertex of $V_1$ for which $K$ of the good $2k$-tuples
containing $a$ are not edges then there
are at least $K - 10^{4k^2} n^{k-1/2}$ bad edges containing $a$.\\
2. If $b$ is a vertex of $V_2$ for which $L$ of the good $2k$-tuples
containing $b$ are not edges then there
are at least $L - (\eta n/5)^{2k-1}$ bad edges containing $b$.
\end{claim}

\noindent {\bf Proof.}\
1. By the preceding remarks $a$ belongs to
at most $\frac{1}{2} {{n-1} \choose {2k-1}} + 10^{4k^2} n^{k-1/2}$
good $2k$-tuples and therefore it belongs to at most
$\frac{1}{2} {{n-1} \choose {2k-1}} + 10^{4k^2} n^{k-1/2} - K$
good edges. Then by equation (\ref{mindeg}) $a$ belongs to at
least $K - 10^{4k^2} n^{k-1/2}$ bad edges.\\
2. From Lemma \ref{kraw}, part (5) $b$ belongs to at most
$\frac{1}{2} {{n-1} \choose {2k-1}} + {{2\cdot\frac{1}{10}\eta n-1} \choose {2k-1}}$
good $2k$-tuples, and the stated bound follows as in (1). \hfill $\Box$

\vspace{0.15cm}
Before proving the next claim we make a remark that will be used on
several occasions without further comment. Suppose $W$ is a bad edge, so that
$|W \cap V_i|$ is even for $i=1,2$. If we partition $W = P \cup Q$ 
with $|P|=|Q|=k$ then $|P \cap V_i|=|Q \cap V_i|$ (mod $2$) for $i=1,2$.
Then for any $k$-tuple $R$ with $|R \cap V_i|=|P \cap V_i|+1$ (mod $2$)
both $2k$-tuples $P \cup R$ and $Q \cup R$ are good. We can obtain
such a $k$-tuple $R \subset V - (P \cup Q)$ by picking any $(k-1)$-tuple,
and then another vertex which, because of parity, is constrained to
lie in some particular $V_i$. This counts each $k$-tuple $k$ times,
so the number of choices for $R$ is at least
$$k^{-1} {{n-2k} \choose {k-1}} \left( \Big( \frac{1}{2} - \frac{1}{10}\eta \Big)n - 3k 
\right)> n^k / (3\cdot k!).$$

\begin{claim}
\label{descend}
Suppose $t\leq k$ and $T$ is a $t$-tuple of vertices
belonging to $\theta n^{2k-t}$ bad edges, for some $\theta > (20k)^k \eta$. Then any
$S \subset T$ with $|S|=t-1$ belongs to at least $(10k)^{-k} \theta n^{2k-t+1}$
good non-edges.
\end{claim}

\noindent {\bf Proof.}\
Write $T = S \cup v$. Consider a bad edge $W$ containing $T$
and a partition $W = T \cup X \cup Y$, where $|X|=k-1$ and
$|Y|=k+1-t$. By the above remark, there are at least $n^k / (3\cdot k!)$ 
$k$-tuples $R$ for
which $v \cup X \cup R$ and $S \cup Y \cup R$ are both good
$2k$-tuples. Note that they can't both be edges, or we would
have a copy of ${\cal C}^{(2k)}_3$. 
Suppose that for at least $\frac{1}{2} \theta n^{2k-t}$ such
$W$ there is a partition $W = T \cup X \cup Y$ for which there
are at least $n^k / 2(3\cdot k!)$ $k$-tuples $R$ for which
$v \cup X \cup R$ is a good non-edge. This clearly gives at least
$\frac{1}{2} \theta n^{k-1}$ choices for $X$. Each such non-edge can
be partitioned in at most ${{2k-1} \choose k}$ ways in the form $v \cup X \cup 
R$, so there are at least
$$ {{2k-1} \choose k}^{-1} \, 
\frac{1}{2} \theta n^{k-1} \, \frac{n^k}{ 2(3\cdot k!)} > (10k)^{-k} \theta 
n^{2k-1}$$
good non-edges containing $v$. Now Claim \ref{goodtobad} shows
that there are at least
$$(10k)^{-k} \theta n^{2k-1} - (\eta n/5)^{2k-1} 
> (20k)^{-k} \theta n^{2k-1} > \eta n^{2k-1}$$
bad edges containing $v$, which contradicts Claim \ref{baddeg}.
It follows that for at least $\frac{1}{2} \theta n^{2k-t}$ such
$W$ and any partition of $W = T \cup X \cup Y$ we have a
good non-edge $S \cup Y \cup R$ for at least $n^k / 2(3\cdot k!)$ $k$-tuples 
$R$.
This gives at least $\frac{1}{2} \theta n^{k+1-t}$ choices for $Y$.
Each such non-edge has at most ${{2k-t+1} \choose k}$ representations
as $S \cup Y \cup R$, so there are at least
$$ {{2k-t+1} \choose k}^{-1} \,
\frac{1}{2} \theta n^{k+1-t} \, \frac{n^k}{ 2(3\cdot k!)} > (10k)^{-k} \theta 
n^{2k-t+1}$$
good non-edges containing $S$. \hfill $\Box$

\vspace{0.15cm}
Suppose for the sake of contradiction that there is some bad edge incident
with $V_1$. Denote the set of bad edges containing some vertex 
$v$ by ${\cal Z}(v)$.
Let $a$ be a vertex in $V_1$ belonging to the maximum number of bad edges
and let $Z=|{\cal Z}(a)|$. Note that $Z>0$.

\begin{claim}
\label{descend2}
Suppose $t \leq k$, $(20k)^k \eta < \phi < (100k)^{-k}$
and ${\cal F}$ is a set of at least $\phi Z n^{-(2k-t)}$ $t$-tuples
containing $a$ such that each $F \in {\cal F}$ is contained in at least $\phi n^{2k-t}$ bad edges.
Then there are at least $\phi^5 Z n^{-(2k-t+1)}$ $(t-1)$-tuples
containing $a$ each of which is contained in at least $\phi^5 n^{2k-t+1}$ bad edges.
\end{claim}

\noindent {\bf Proof.}\
Let ${\cal G}$ be the set of $(t-1)$-tuples containing $a$ that are contained in
a member of ${\cal F}$. By claim \ref{descend} each $G \in {\cal G}$ is contained
in at least $(10k)^{-k} \phi n^{2k-t+1}$ good non-edges. Each such good non-edge
is counted by at most ${{2k-1} \choose {t-2}}$ different $G$'s, so there are at least
${{2k-1} \choose {t-2}}^{-1} |{\cal G}| (10k)^{-k} \phi n^{2k-t+1}
> (40k)^{-k} |{\cal G}| \phi n^{2k-t+1}$ good non-edges containing $a$.
Since $a \in V_1$ Claim \ref{goodtobad} gives at least
$$(40k)^{-k} |{\cal G}| \phi n^{2k-t+1} - 10^{4k^2} n^{k-1/2}
> (50k)^{-k} |{\cal G}| \phi n^{2k-t+1}$$
bad edges containing $a$, so by definition of $Z$ we get 
$|{\cal G}|  < (50k)^k \phi^{-1} Z n^{-(2k-t+1)}$.
Let ${\cal G}_1 \subset {\cal G}$ consist of those $G$ that belong to at
least $\phi^3 n$ members of ${\cal F}$. Then 
$$\phi Z n^{-(2k-t)} \leq |{\cal F}| < |{\cal G}_1| n + |{\cal G}| \phi^3 n
< |{\cal G}_1| n + (50k)^k \phi^2 Z n^{-(2k-t)}$$
so $|{\cal G}_1| > \phi^5 Z n^{-(2k-t-1)}$ with room to spare. For each $G \in 
{\cal G}_1$
there are at least $\phi^3 n$ sets of ${\cal F}$ each contributing
$\phi n^{2k-t}$ bad edges containing $G$. Each such bad edge is counted by at most
$2k-t+1$ different $F \in {\cal F}$, so $G$ 
belongs to at least $(2k-t+1)^{-1} \phi^3 n \cdot \phi n^{2k-t} > \phi^5 n^{2k-t+1}$ bad edges. \hfill $\Box$

\vspace{0.15cm}
Let ${\cal Z}_1(a)$ be those bad edges $W$ containing $a$ 
for which there is some partition into two $k$-tuples $W = P \cup Q$ with $a 
\in P$
so that there are at least $n^k / 2(3\cdot k!)$ $k$-tuples $R$ for which
$P \cup R$ is a good non-edge. Let ${\cal Z}_2(a)={\cal Z}(a)-{\cal Z}_1(a)$, and write 
$Z_i = |{\cal Z}_i(a)|$ for $i=1,2$. Then one of $Z_1$,$Z_2$ is at least $Z/2$.

{\bf Case 1:}\,
Suppose $Z_1 \geq Z/2$. Let ${\cal P}$ be the (non-empty) set of $k$-tuples $P$
containing $a$ such that there is some edge $P \cup Q$ in ${\cal Z}_1(a)$,
and $P \cup R$ is a good non-edge for at least  $n^k / 2(3k!)$ $k$-tuples $R$.
Each such good non-edge is counted by at most ${{2k-1} \choose {k-1}}$
different $P$'s, so there are at least
${{2k-1} \choose {k-1}}^{-1} |{\cal P}| n^k / 2(3\cdot k!) > (10k)^{-k} |{\cal 
P}| n^k$ 
good non-edges containing $a$. Now Claim \ref{goodtobad} gives at least
$(10k)^{-k} |{\cal P}| n^k - 10^{4k^2} n^{k-1/2} > (20k)^{-k} |{\cal P}| n^k$
bad edges containing $a$, so by definition of $Z$, $|{\cal P}|  < (20k)^k  Z 
n^{-k}$. On the other hand, 
let ${\cal P}_1 \subset {\cal P}$ consist of those $P$ that belong to at
least $\frac{1}{10} (20k)^{-k} n^k$ bad edges. Then
$$ Z/2 \leq Z_1 < |{\cal P}_1| n^k + |{\cal P}| \frac{1}{10} (20k)^{-k} n^k
< |{\cal P}_1| n^k + Z/10$$
so $|{\cal P}_1| > 0.4\ Z n^{-k}$. Now apply
Claim \ref{descend2} $k-1$ times, starting with $t=k$ and
$\phi = (100k)^{-k}$. We deduce that $a$ belongs to
at least $\phi^{5^{k-1}} n^{2k-1} > \eta n^{2k-1}$ bad edges, which contradicts
Claim \ref{baddeg}.

{\bf Case 2:}\,
Now suppose $Z_2 \geq Z/2$. Note that every bad edge containing $a$
contains at least one other point of $V_1$, so there is
some $b \in V_1$ belonging to at least $Z_2/n$
edges of ${\cal Z}_2(a)$. Fix one such $b$. Let 
${\cal X}$ be the set of $(k-1)$-tuples $X$ for which there
exists a $(k-1)$-tuple $Y$ such that $W=a \cup b \cup X \cup Y$
is an edge of ${\cal Z}_2(a)$. By definition of ${\cal Z}_2(a)$ for any 
such partition of $W$, there
are at least $n^k / 2(3\cdot k!)$ $k$-tuples $R$ such that $b \cup X \cup R$
is a good non-edge. This gives at least
$\frac{n^k}{2(3\cdot k!)}|{\cal X}|>(10k)^{-k} |{\cal X}| n^k$ good non-edges 
containing $b$, and
since $b \in V_1$ Claim \ref{goodtobad} gives at least 
$(10k)^{-k} |{\cal X}| n^k-10^{4k^2}n^{k-1/2}>(20k)^{-k} |{\cal X}| n^k$ bad
edges containing $b$. Thus, by definition of $Z$, $|{\cal X}| < (20k)^k Z 
n^{-k}$. Note that
each edge in ${\cal Z}_2(a)$ that contains $b$ is obtained by
picking a pair of $(k-1)$-tuples in ${\cal X}$, so
$Z/(2n) \leq Z_2/n \leq {{|{\cal X}|} \choose 2} < \frac{1}{2} (20k)^{2k} Z^2 n^{-2k}$.
Therefore $Z > (20k)^{-2k} n^{2k-1} > \eta n^{2k-1}$, which contradicts
Claim \ref{baddeg}.

We conclude that there are no bad edges incident to the vertices of $V_1$, i.e. all bad edges
are entirely contained in $V_2$. As in the case $k=2$ this gives a more
precise bound on $t$, defined by $|V_1|=n/2+t$, $|V_2|=n/2-t$. If
$|t| \geq 20^k \sqrt{n}$ then Lemma \ref{kraw}, part (4) gives
$d_{2k}(n,t) < \frac{1}{2} {{n-1} \choose {2k-1}} - 20^k n^{k-1/2}$. 
This is a contradiction, since the
vertices of $V_1$ only belong to good edges, of which there
are at most $d_{2k}(n,t)<\delta(H)$. Therefore $|t|<20^k \sqrt{n}$.
Now Lemma \ref{kraw}, part (2) gives
$$d_{2k}(n,-t) < \frac{1}{2} {{n-1} \choose {2k-1}} +
\big(10(20^k)^2\big)^k n^{k-1/2}<
\frac{1}{2} {{n-1} \choose {2k-1}} + 10^{4k^2} n^{k-1/2}.$$
As we remarked earlier, this bound allows us to repeat the above argument
interchanging $V_1$ and $V_2$, so we deduce that there are no bad edges
incident with $V_2$ either, i.e. all edges are good. Then by definition of $b_{2k}(n)$
we have $e(H) \leq b_{2k}(n)$, with equality only when $H$ is a ${\cal B}^{(2k)}_n$,
so the theorem is proved. \hfill $\Box$

\section{Hypergraphs without ${\cal C}^{(4)}_r$}

We recall that ${\cal C}^{(2k)}_r$ is the $2k$-uniform hypergraph obtained 
by letting $P_1,\cdots,P_r$ be pairwise disjoint sets of size $k$ and 
taking as edges all sets $P_i \cup P_j$ with $i \neq j$. In this section
we will be concerned with the case $k=2$ and general $r$.

Sidorenko \cite{Si} showed that the Tur\'an density of ${\cal C}^{(2k)}_r$ 
is at most $\frac{r-2}{r-1}$. This is a consequence of Tur\'an's theorem applied
to an auxiliary graph $G$ constructed from a $2k$-uniform hypergraph $H$ 
of order $n$.
The vertices of $G$ are the $k$-tuples of vertices of $H$, and
two $k$-tuples $P_1$,$P_2$ are adjacent if $P_1 \cup P_2$ is an edge of $H$.
It is easy to see that the graph $G$ has ${n \choose k}$ vertices,
$\frac{1}{2}{2k \choose k}e(H)$ edges and contains no $K_r$. Thus the
upper bound on the number of edges of $H$ follows immediately from
Tur\'an's theorem. The following construction from \cite{Si}
gives a matching lower bound when $r$ is of the form $2^p+1$.

Let $W$ be a vector space of dimension $p$ over the field 
$GF(2)$, i.e. the finite field with $2$
elements $\{0,1\}$. Partition a set of vertices $V$ as $\bigcup_{w \in W} 
V_w$, $|V_w| =|V|/(r-1)$.
Given $t$ and a $t$-tuple of vertices $X = x_1 \cdots x_t$ with $x_i \in V_{w_i}$ 
we define $\Sigma X = \sum_1^t w_i$. Define a $2k$-uniform hypergraph $H$, where
a $2k$-tuple $X$ is an edge iff $\Sigma X \neq 0$. Observe that this
doesn't contain a copy of ${\cal C}^{(2k)}_r$. Indeed, if 
$P_1,\cdots,P_r$ are disjoint $k$-tuples
then there is some $i \neq j$ with $\Sigma P_i = \Sigma P_j$ (by the
pigeonhole principle). Then $\Sigma (P_i \cup P_j) =\Sigma P_i+ \Sigma P_j
= 0$, so $P_i \cup P_j$ is
not an edge. 

This construction depends essentially on an algebraic structure, which
only exists for certain values of $r$. Perhaps surprisingly, we will
show that this is an intrinsic feature of the problem, by proving
Theorem \ref{c4r}, which gives a stronger upper bound on the Tur\'an density of 
${\cal C}^{(4)}_r$, when
$r$ is not of the form $2^p+1$. We make no attempt to optimize the 
constant in this bound.

In addition, our proof of this theorem implies that, 
for $r=2^p+1$, any ${\cal C}^{(4)}_r$-free $4$-uniform hypergraph
with density $\frac{r-2}{r-1}- o(1)$ looks approximately like 
Sidorenko's construction.

\begin{coro}
\label{c42}
Let $r=2^p+1$ be an integer and let $W$
be a $p$-dimensional vector space over the field $GF(2)$.
For every $\epsilon>0$ there is $\eta>0$ so that if $H$ is a
${\cal C}^{(4)}_r$-free
$4$-uniform hypergraph with $e(H)>\frac{r-2}{r-1}{n \choose 4}-\eta n^4$ 
then there is a partition of the vertex set as $\bigcup_{w \in W} V_w$
such that all but $\epsilon n^4$ edges $X$ of $H$ satisfy $\Sigma X\not 
=0$.
\end{coro}

The rest of this section is organized as follows.
In the first subsection we will prove a lemma 
showing that  certain
edge-colourings of the complete graph $K_s$ exist only if $s$ is 
a power of $2$.
In the following subsection we will recall a proof of the
Simonovits stability theorem so that we can calculate some explicit 
constants. The final subsection contains the proof
of Theorem \ref{c4r}. 

\subsection{A lemma on edge-colourings of a complete graph}

\begin{lemma}
\label{algebra}
Suppose that we have a colouring of the edges of the complete graph $K_s$
in $s-1$ colours, so that every colour is a matching and each
subset of $4$ vertices spans edges of either $3$ or $6$ different colours.
Then $s=2^p$ for some integer $p$.
\end{lemma}

\noindent {\bf Proof.}\
Since the number of colours is $s-1$, every colour is a matching and the 
total number of edges in $K_s$ is $s(s-1)/2$ it is easy to see that
every colour is a perfect matching. Also, if $wx$ and $yz$ are
disjoint edges of the same colour, then by hypothesis only
$3$ different colours appear on $wxyz$, so $wy$ and $xz$ have
the same colour, as do $xy$ and $wz$. Denote the set of colours by
$C = \{c_1,\cdots,c_{s-1}\}$. We define
a binary operation $+$ on $C$ using the following rule. Pick a
vertex $x$. Given $c_i$ and $c_j$ let $e_i=xy_i$ and $e_j=xy_j$ be
the edges incident with $x$ with these colours. These edges
exist, as each colour is a perfect matching. Define $c_i + c_j$ to
be the colour of $y_i y_j$.

To see that this is well-defined, let $x'$ be another vertex
and suppose $e_i'=x'y_i'$ has colour $c_i$ and 
$e_j'=x'y_j'$ has colour $c_j$. If $y_i=y'_j$ then opposite
edges of $xy_jy_ix'$ have the same colours, so $x'y_j$ has
colour $c_i$, i.e. $y_j=y'_i$ and there is nothing to prove.
Therefore we can assume that all $y_i, y_j, y'_i, y'_j$ are distinct. 
Consider the $4$-tuple $xx'y_iy_i'$.
Since $xy_i$ and $xy_i'$ have the same colour 
we deduce that $xx'$ and $y_i y_i'$ have
the same colour. Similarly $xx'$ and $y_j y_j'$ have
the same colour, from which we see that $y_i y_i'$ and $y_j y_j'$ have
the same colour. Now looking at $y_i y_i' y_j y_j'$ we see that
$y_i y_j$ and $y_i' y_j'$ have the same colour, so
$c_i + c_j$ is well-defined.

Let $D$ be a set obtained by adjoining another element called
$\bf 0$ to $C$. Extend $+$ to an operation on $D$ by defining
${\bf 0}+d=d+{\bf 0}=d$ and $d+d={\bf 0}$ for all $d \in D$. We claim that 
$(D,+)$
is an abelian group. Note that $+$ is commutative by definition, ${\bf 0}$ 
is
an identity and inverses exist. It remains to show
associativity, i.e. for any $d_1,d_2,d_3$ we have
$(d_1+d_2)+d_3=d_1+(d_2+d_3)$. This is immediate if any
of the $d_i$ are ${\bf 0}$ or if they are all equal. If $d_1=d_2 \neq d_3$
then $d_1+d_2={\bf 0}$ and there is a triangle with colours 
$d_1,d_3,d_1+d_3$,
so $d_1+(d_2+d_3)=d_3$ as required. The same argument
applies when $d_2=d_3 \neq d_1$. If $d_1=d_3$ then $d_1+d_2=d_2+d_3$
by commutativity, and so $(d_1+d_2)+d_3=d_1+(d_2+d_3)$ also
by commutativity. So we can assume that the $d_i$ are pairwise
distinct and non-zero. Pick a vertex $x$, let $xy_1$ be
the edge of colour $d_1$ and $xy_2$ the edge of colour $d_2$.
Let $y_2 z$ have colour $d_3$. We can suppose $z \neq y_1$,
otherwise $d_1+d_2=d_3$ and $d_2+d_3=d_1$ and 
$(d_1+d_2)+d_3=d_1+(d_2+d_3)={\bf 0}$. Now $y_1 y_2$ has colour
$d_1 + d_2$ and $xz$ has colour $d_2 + d_3$. Consider
the edge $y_1 z$. From the triangle it forms with $x$ we see
that it has colour $d_1+(d_2+d_3)$ and from the triangle
with $y_2$ we see that it has colour $(d_1+d_2)+d_3$. This
proves associativity, so $D$ is an abelian group.

Finally, note that every non-zero element has order $2$, so
$D$ is in fact a vector space over the field with $2$ elements.
If $p$ is its dimension then $s=|D|=2^p$. \hfill $\Box$

\subsection{The Simonovits stability theorem}
In this subsection we will recall a proof of the
Simonovits stability theorem \cite{S1} so that we can calculate some
explicit constants.
Let $T_s(N)$ be the $s$-partite {\em Tur\'an graph} on $N$
vertices, i.e. a complete $s$-partite graph with part sizes
as equal as possible. Write $t_s(N)$ for the number of edges
in $T_s(N)$. Then Tur\'an's theorem states that any
$K_{s+1}$-free graph on $N$ vertices has at most $t_s(N)$ edges,
with equality only for $T_s(N)$. It is easy to show that
$\frac{s-1}{s}N^2/2 - s < t_s(N) \leq \frac{s-1}{s}N^2/2$.

\begin{prop}
\label{stab1}
Suppose $G$ is a $K_{s+1}$-free graph on $N$ vertices with
minimum degree $\delta(G) \geq \big( 1 - \frac{1}{s} - \alpha \big)N$ and 
$\alpha<1/s^2$.
Then there is a partition of the vertex set of $G$ as
$V(G) = U_1 \cup \cdots U_s$ with $\sum e(U_i) < s\alpha N^2$.
\end{prop}

\noindent {\bf Proof.}\
By Tur\'an's theorem $G$ contains a copy of $K_s$; let
$A=\{a_1,\cdots,a_s\}$ be its vertex set. Note that any vertex $x$ not in $A$
has at most $s-1$ neighbours in $A$, or we get a $K_{s+1}$.
Let $B$ be those vertices with exactly $s-1$ neighbours in $A$,
and $C=V(G)-A-B$. Partition $A \cup B$ as $U_1 \cup \cdots \cup U_s$ where
$U_i$ consists of those vertices adjacent to $A-a_i$. Then there are no edges
inside any $U_i$, as if $xy$ is such an edge then $xy+A-a_i$ forms
a $K_{s+1}$. Distribute the vertices of $C$
arbitrarily among the $U_i$. Counting edges between $A$ and $V-A$ gives
$$s\delta(G) \leq e(A,V-A) \leq (s-1)|B|+(s-2)|C| = (s-1)(N-s)-|C|$$
so $|C| \leq s\alpha N - s(s-1)$. Therefore $\sum e(U_i) < s\alpha N^2$.
\hfill $\Box$
 
\begin{theo}
\label{stab2}
Suppose $G$ is a $K_{s+1}$-free graph on $N$ vertices with
at least $\Big( \frac{s-1}{2s} - c \Big) N^2$ edges and $c<1/(4s^4)$. Then
there is a partition of the vertex set of $G$ as
$V(G) = U_1 \cup \cdots \cup U_s$ with $\sum e(U_i) < (2s+1)\sqrt{c}\ N^2$.
\end{theo}

\noindent {\bf Proof.}\
Construct a sequence of graphs $G=G_N,G_{N-1},\cdots$
where if $G_m$ has a vertex of degree at most
$\Big( 1 - \frac{1}{s} - 2\sqrt{c} \Big)m$ then we delete it
to get $G_{m-1}$. Suppose can delete $\sqrt{c}\ N$ vertices by
this process and reach a graph $G_{(1-\sqrt{c})N}$. Then
$G_{(1-\sqrt{c})N}$ is $K_{s+1}$-free and has at least
$$\left( \frac{s-1}{2s} - c - \sqrt{c} \Big( 1 - \frac{1}{s} - 2\sqrt{c} \Big) \right) N^2
> \frac{s-1}{2s}(1-\sqrt{c})^2N^2$$
edges. This contradicts Tur\'an's theorem, so 
the sequence terminates at some $G_m$ with $m \geq (1-\sqrt{c})N$
and minimum degree at least $\Big( 1 - \frac{1}{s} - 2\sqrt{c} \Big)m$.
By Proposition \ref{stab1} there is a partition
$V(G_m) = U_1 \cup \cdots \cup U_s$ with $\sum e(U_i) < 2s\sqrt{c}\ N^2$.
Now distribute the $\sqrt{c} N$ deleted vertices arbitrarily among
the $U_i$. Then $\sum e(U_i) < (2s+1)\sqrt{c}\ N^2$. \hfill $\Box$

\subsection{Proof of Theorem \ref{c4r}}

Let $V$ be the vertex set of $H$. Define a graph $G$ whose vertices are all pairs 
in $V$, where the pairs $ab$ and $cd$ are adjacent exactly when $abcd$ is an
edge of $H$. Since $H$ is ${\cal C}^{(4)}_r$-free we see that $G$ is 
$K_r$-free.
Also, each edge of $H$ creates exactly $3$ edges in $G$ (corresponding
to the $3$ ways of breaking a $4$-tuple into pairs) so
$$e(G) > 3\Big(\frac{r-2}{r-1} - 10^{-33}r^{-70} \Big) {n \choose 4}
> \Big( \frac{r-2}{2(r-1)} - 10^{-33}r^{-70} \Big) N^2,$$
where $N={n \choose 2}$.

Applying Theorem \ref{stab2} with $s=r-1$ gives a partition of the 
pairs of vertices in $V$
as $\bigcup_1^{r-1} P_i$ with $\sum_1^{r-1} e(P_i) < 10^{-16}r^{-34} 
N^2$. If there is some $P_i$ with
$|P_i| < \big(\frac{1}{r-1} - 10^{-3}r^{-7}\big)N$ then
\begin{eqnarray*}
\frac{e(G)}{N^2} &<& \frac{{r-2 \choose 2}}{(r-2)^2} \left( 
\frac{r-2}{r-1} + 
10^{-3}r^{-7}\right)^2
+ \left(\frac{1}{r-1} - 10^{-3}r^{-7}\right)\left(\frac{r-2}{r-1} + 
10^{-3}r^{-7}\right)
+ 10^{-16}r^{-34}\\
&<& \frac{r-2}{2(r-1)} - 10^{-6}r^{-14}/2 + 10^{-16}r^{-34}\,.
\end{eqnarray*}
This is a contradiction so $|P_i| \geq \big(\frac{1}{r-1} - 
10^{-3}r^{-7}\big)N$ for all $i$. Also if some 
$|P_i| > \big(\frac{1}{r-1} + 10^{-3}r^{-6}\big)N$, then there is $j$ 
such that $|P_j|<\big(\frac{1}{r-1} -10^{-3}r^{-7}\big)N$. Therefore 
for all $i$
\begin{equation}
\label{sizeP_i}
\left||P_i|-\frac{1}{r-1}N\right|\leq 10^{-3}r^{-6}n^2\,.
\end{equation}

Note that all but at most $10^{-16}r^{-34} n^4$ edges of $H$ are formed 
by taking
a pair from $P_i$ and a pair from $P_j$ with $i \neq j$.
We think of the $P_i$ as a colouring of pairs.
A $4$-tuple $abcd$ will be called {\em properly coloured} if either 

\noindent
(i) $abcd$ is an edge and each of the $3$ sets  
$\{ab,cd\}$,$\{ac,bd\}$,$\{ad,bc\}$
contains two pairs with different colours, or

\noindent
(ii) $abcd$ is not an edge and each of the $3$ sets 
$\{ab,cd\}$,$\{ac,bd\}$,$\{ad,bc\}$
consists of two pairs with the same colour.

An improperly coloured $4$-tuple is either an edge that is the union of 
two pairs of the same colour or
a non-edge which is the union of two pairs with different colours. 
There are at most $10^{-16}r^{-34}N^2$ of the former $4$-tuples, and the 
number of latter is at most
$$\frac{r-2}{2(r-1)}N^2-\Big(e(G)-10^{-16}r^{-34}N^2\Big)\leq 
\Big(10^{-16}r^{-34} + 
10^{-33}r^{-70}\Big) 
n^4\,.$$
Therefore
all but $10^{-15}r^{-34} n^4$ $4$-tuples are properly coloured. Call a 
pair $ab$ {\em bad} if there are at least $10^{-12}r^{-32} n^2$ pairs $cd$ 
such that $abcd$
is improperly coloured; otherwise call it {\em good}.
Then there are at most 
${4 \choose 2}(10^{-15}r^{-34} n^4)/(10^{-12}r^{-32} n^2)<10^{-2}r^{-2} n^2$ 
bad pairs.

Consider a graph on $V$ whose edges are the pairs in $P_1$. As noted 
in (\ref{sizeP_i}) it has at least $\frac{1}{r-1}N - 10^{-3}r^{-6} n^2$ edges.
For vertices $a$ and $b$ in $V$, 
let $d(a)$ denote the degree of $a$
and $d(a,b)$ the {\em codegree} of $a$ and $b$ 
(i.e. the size of their common neighbourhood.) Then
$$\sum_{a,b \in V} d(a,b) = \sum_{c \in V} {{d(c)} \choose 2}
\geq n {{\sum d(c)/n} \choose 2}=
n { {2|P_1|/n} \choose 2}  > \frac{1}{5r^2} nN.$$
Suppose there are at most $m$ pairs $(a,b)$ for which 
$d(a,b) > \frac{n}{10r^2}$. Then 
$\frac{1}{5r^2} nN < \sum d(a,b) \leq mn + N\frac{n}{10r^2}$,
so $\frac{1}{10r^2} N < m$, i.e. there are at least
$\frac{1}{10r^2} {n \choose 2}$ pairs $(a,b)$ for
which $d(a,b) > \frac{n}{10r^2}$. At least one such pair is good,
as the number of bad pairs is at most $10^{-2}r^{-2} n^2 < \frac{1}{20r^2} n^2$.
Let $(a,b)$ be such a pair and suppose it belongs to $P_t$.

Let $B$ be the set of pairs $cd$ for which $abcd$ is improperly
coloured. Since $ab$ is good we have $|B| \leq 10^{-12}r^{-32} n^2$.
Therefore there are at most $|B|{n \choose 2} < 10^{-12}r^{-32} n^4$
$4$-tuples of vertices that contain any pair of $B$. We will call
a $4$-tuple {\em normal} if it is properly coloured and
does not contain a pair from $B$; otherwise we call it {\em abnormal}. Then all but
at most $10^{-12}r^{-32} n^4+10^{-15}r^{-34} n^4<10^{-11}r^{-32} n^4$ 
$4$-tuples are normal.

Partition the vertices of $V-ab$ into $(r-1)^2$ sets $U_{ij}$, where
$c$ is in $U_{ij}$ iff $ac \in P_i$ and $bc \in P_j$. Then
by the above discussion $|U_{11}| \geq \frac{n}{10r^2}$.
Now we claim that for $i \neq j$ we have $|U_{ij}| < 10^{-3}r^{-11} n$.
For suppose that $|U_{ij}| \geq 10^{-3}r^{-11} n$. Let $P_k$ be the colour 
that appears most frequently among pairs joining vertices of $U_{11}$
to $U_{ij}$. Then there are at least $\frac{1}{r-1}|U_{11}||U_{ij}|$ 
pairs of color $P_k$ with one endpoint in $U_{11}$ and the other in 
$U_{ij}$. 
Consider the $4$-tuples of the form $c_1 c_2 d_1 d_2$, with 
$c_1,c_2 \in U_{11}$, $d_1,d_2 \in U_{ij}$ and $c_1d_1,c_2d_2 \in P_k$.
There are at least 
$$\left( \frac{1}{r-1} |U_{11}||U_{ij}| \right)
\left( \frac{1}{r-1} |U_{11}||U_{ij}|-2n \right)/4>
r^{-2}\Big(10^{-1}r^{-2}n \cdot 10^{-3}r^{-11} n\Big)^2/4
> 10^{-11}r^{-32} n^4$$
such $4$-tuples, so some $c_1 c_2 d_1 d_2$ is normal.
By definition of normality each of its pairs forms a properly
coloured $4$-tuple with $ab$. Since $ac_1$ and $bc_2$ are in $P_1$ and
$ab$ is in $P_t$ we deduce that $c_1 c_2$ is in $P_t$ as well. 
Also $ad_1 \in P_i$, $bd_2 \in P_j$ and $i \neq j$, so
$d_1 d_2$ cannot be in $P_t$. But $c_1d_1$ and $c_2d_2$ both belong to $P_k$
so $c_1 c_2 d_1 d_2$ is improperly coloured. This contradicts
the definition of normality, so we do have $|U_{ij}| < 10^{-3}r^{-11} n$.

For convenience write $U_i = U_{ii}$. Then
all but at most $(r-1)^210^{-3}r^{-11} n\leq 10^{-3}r^{-9} n$ vertices 
belong to one of the $U_i$.
Suppose $cd$ is a pair such that $abcd$ is properly coloured.
Since $ab$ is a good pair, this is the case for all but at most 
$10^{-11}r^{-32} n^2$ pairs $cd$.
If $c$ and $d$ both belong to some $U_i$ then $ac$ and $bd$ both
have colour $i$. Since $ab$ has colour $t$ we see that $cd$ has colour $t$.
Similarly, if $c \in U_i$ and $d \in U_j$ with $i \neq j$ we
see that $cd$ cannot have colour $t$.

Let $E_i$ denote the pairs with both endpoints in $U_i$,
so that $|E_i| = {{|U_i|} \choose 2}$. By the above discussion, all 
but at most
$10^{-12}r^{-32} n^2$ pairs in $\cup_i E_i$ belong to $P_t$. Suppose
$|U_i| < \big( \frac{1}{r-1} - 10^{-1}r^{-3} \big)n$ for some $i$,
so that
\begin{eqnarray*} 
\sum |E_i| &>& \left( \frac{1}{r-1} - 10^{-1}r^{-3} 
\right)^2\frac{n^2}{2}+
(r-2)\left(\frac{1-1/(r-1)+10^{-1}r^{-3} - 10^{-3}r^{-9}}{r-2}
\right)^2\frac{n^2}{2}-O(n)\\
&>& \left(\frac{1}{r-1} + 10^{-2}r^{-6}\right)\frac{n^2}{2} - O(n).
\end{eqnarray*}
By (\ref{sizeP_i}), this gives the following contradiction.
$$\frac{1}{r-1}{n \choose 2} + 10^{-3}r^{-6} n^2
\geq |P_t| \geq \sum_i |E_i| - 10^{-12}r^{-32} n^2>
\frac{1}{r-1}\frac{n^2}{2}+\frac{10^{-2}r^{-6}}{3}n^2.$$
Therefore $|U_i|\geq \big( \frac{1}{r-1} - 10^{-1}r^{-3} \big)n$
for each $i$. 

Let $E_{ij}$ denote the edges with one endpoint in $U_i$
and the other in $U_j$. We claim that one colour is dominant
among these edges, i.e. there is some $q$ such that all but $10^{-2}r^{-4} n^2$ 
edges of $E_{ij}$ belong to $P_q$. Indeed, suppose that there are colours
$q_1$ and $q_2$ for which there are at least $10^{-2}r^{-4} n^2$ edges in 
$E_{ij}$ of color $q_i$ for $i=1,2$. Then there are at least 
$(10^{-2}r^{-4} n^2)(10^{-2}r^{-4} n^2 - 2n) > 10^{-11}r^{-32} n^4$ 
$4$-tuples $c_1 c_2 d_1 d_2$
with $c_1,c_2$ in $U_i$, $d_1,d_2$ in $U_j$ and $c_i d_i$ of colour $q_i$.
At least one such $4$-tuple $c_1 c_2 d_1 d_2$ is normal, since there at most
$10^{-11}r^{-32} n^4$ abnormal $4$-tuples. But then $c_1 c_2$ and $d_1 d_2$
both have colour $t$, so $c_1 c_2 d_1 d_2$ is improperly coloured,
which is a contradiction.

Consider the complete graph $K_{r-1}$ on the vertex set $\{1,\cdots,r-1\}$ 
and colour edge $ij$ with the dominant colour of $E_{ij}$. We show
that this colouring satisfies the hypotheses of Lemma \ref{algebra}.
First of all we show that colour $t$ doesn't occur in this edge-coloring 
of $K_{r-1}$, i.e. there
are only $r-2$ colours. Suppose $ij$ has colour $t$. Then

\begin{eqnarray*}
\frac{1}{r-1}{n \choose 2} + 10^{-3}r^{-6} n^2& \geq&
|P_t| \geq \sum_i |E_i| - 10^{-12}r^{-32} n^2 + |U_i||U_j| - 10^{-2}r^{-4} n^2\\
&\geq& (r-1) { {( \frac{1}{r-1} - 10^{-1}r^{-3} ) n } \choose 2}
+ \left(\Big( \frac{1}{r-1} - 10^{-1}r^{-3} \Big) n\right )^2 - 
10^{-1}r^{-4} n^2\\
&\geq& \frac{1}{r-1}\,\frac{n^2}{2}+\frac{n^2}{(r-1)^2}-
\frac{r+2}{r-1}10^{-1}r^{-3}n^2>\frac{1}{r-1}\,
\frac{n^2}{2}+\frac{n^2}{r^2},
\end{eqnarray*}
is a contradiction. 

Now suppose that some colour $\ell$ is not a matching, i.e.
there are edges $ij$ and $ik$ in $K_{r-1}$ both of colour $\ell$. Then
all but at most $2\cdot 10^{-2}r^{-4} n^2$ pairs of $E_{ij} \cup E_{ik}$
have colour $\ell$. Consider the $4$-tuples of the form $c_1 c_2 d e$, with 
$c_1,c_2 \in U_i$, $d \in U_j$ and $e \in U_k$, such that
$c_1d,c_2d,c_1e,c_2e$ all have colour $\ell$. There
are at least 
$${{|U_i| \choose 2}}|U_j||U_k| - 
2\cdot 10^{-2}r^{-4} n^2{n \choose 2} \geq\frac{1}{2}
\left(\frac{1}{r-1} - 10^{-1}r^{-3}\right)^4n^4-O(n^3)-
10^{-2}r^{-4} n^4 
 > 10^{-11}r^{-32} n^4$$
such $4$-tuples, so there is one such $c_1 c_2 d e$ which is normal.
But then $c_1 c_2$ has colour $t$ and $de$ cannot have
colour $t$, since by normality $abde$ is properly colored.
Therefore $c_1 c_2 d e$  is improperly coloured. This is
a contradiction, so each colour forms a matching.

It remains to show that if some $4$ vertices $x_1 x_2 x_3 x_4$ in 
$K_{r-1}$ do not span $6$ different colours then they span only $3$
colours. Suppose that $x_1 x_2$ and $x_3 x_4$ have colour $\alpha$,
$x_1 x_3$ has colour $\beta$ and $x_2 x_4$ has colour $\gamma$.
Recall that all but at most $10^{-2}r^{-4} n^2$ pairs in $E_{x_ix_j}$ 
have the corresponding color of $x_ix_j$.
Consider the $4$-tuples in $H$ of the form $c_1 c_2 c_3 c_4$ with
$c_i \in U_{x_i}$ such that $c_i c_j$ has the same colour as $x_i x_j$.
There are at least 
$$\prod_1^4 |U_{x_i}| - 4\cdot 10^{-2}r^{-4} n^2{n \choose 2}>
\left(\frac{1}{r-1} - 10^{-1}r^{-3}\right)^4n^4-
2\cdot 10^{-2}r^{-4} n^4
 > 10^{-11}r^{-32} n^4$$
such $4$-tuples. Since the the number of abnormal
$4$-tuples is at most $10^{-11}r^{-32} n^4$,  some such $c_1 c_2 c_3 c_4$ should be  normal. 
Then $\beta = \gamma$, or $c_1 c_2 c_3 c_4$ would be improperly 
coloured. We
see that opposite edges of $x_1 x_2 x_3 x_4$ have the same
colour. Therefore we can apply Lemma \ref{algebra} with
$s=r-1$ to deduce that $r-1$ is of the form $2^p$. \hfill $\Box$

Finally it is not difficult to check that, when $r=2^p+1$, the above 
arguments together with the proof of Lemma \ref{algebra} imply Corollary 
\ref{c42}.

\section{Concluding remarks}

Among the various techniques that we used in this paper, the stability
approach stands out as one that should be widely applicable in
extremal combinatorics. The process of separating the argument into a stability
stage and a refinement stage focuses attention on the
particular difficulties of each, and often leads to progress where the raw problem
has appeared intractable. For recent examples we refer to our proofs of the 
conjecture of S\'os on the Tur\'an number of the Fano plane \cite{KS},
and a conjecture of
Yuster on edge colorings with no monochromatic cliques \cite{ABKS}.

Our methods probably apply to ${\cal C}^{(2k)}_r$ for
general $k$ when $r$ is of the form $2^p+1$, although
the reader who has grappled with the thornier aspects of this paper will
note the formidable technical difficulties that would arise. 
It would be far more interesting to say more about the behaviour
of the Tur\'an density of ${\cal C}^{(2k)}_r$
for general $r$. Even ${\cal C}^{(4)}_4$ presents an enigma for which
there is no obvious plausible conjecture. We find it remarkable that the 
seemingly similar hypergraphs ${\cal C}^{(4)}_3$ and ${\cal C}^{(4)}_5$ are 
actually distinguished 
from ${\cal C}^{(4)}_4$ by a hidden algebraic feature, so are loathe even
to speculate on the nature of the best construction for this case.


\begin{thebibliography}{99}

\bibitem{ABKS} N. Alon, J. Balogh, P. Keevash and B. Sudakov,
The number of edge colorings with no monochromatic cliques, submitted.

\bibitem{Fr} P. Frankl, 
Asymptotic solution of a Tur\'an-type problem. 
{\em Graphs and Combinatorics} 6 (1990), 223--227. 

\bibitem{F1} Z. F\"uredi,
Tur\'an type problems, in: {\em  Surveys in combinatorics}, 
London Math. Soc. Lecture Note Ser. 166, Cambridge Univ. Press, 
Cambridge, 1991, 253--300

\bibitem{KS} P. Keevash and B. Sudakov, 
The exact Tur\'an number of the Fano plane, submitted.

\bibitem{KL} I. Krasikov and  S. Litsyn, 
Survey of binary Krawtchouk polynomials,
in: {\em  Codes and association schemes}, DIMACS Ser. Discrete Math.
Theoret. Comput. Sci., 56, Amer. Math. Soc., Providence, RI, 2001,
199--211.

\bibitem{Le}
V. Levenshtein, Krawtchouk polynomials and universal bounds 
for codes and designs in Hamming spaces, {\em IEEE Trans. Inform. Theory} 
41 (1995), 1303--1321.

\bibitem{Lo} L. Lov\'asz,
{\bf Combinatorial Problems and Exercises},
North-Holland, Amsterdam, 1993.

\bibitem{Si} A. Sidorenko,
An analytic approach to extremal problems for graphs and hypergraphs, in:
{\em Extremal problems for finite sets}, Bolyai 
Soc. Math. Stud. 3, Budapest, 1994, 423--455.

\bibitem{Si1} A. Sidorenko, 
What we know and what we do not know about Tur\'an numbers,
{\em Graphs and Combinatorics} 11 (1995), 179--199.

\bibitem{S1} M.\ Simonovits, A method for solving extremal problems in
graph theory, stability problems, in: {\em Theory of Graphs (Proc.\ 
Colloq.\ Tihany, 1966)\/},  Academic Press,
New York, and Akad.\ Kiad\'o, Budapest, 1968, 279--319.

\end{thebibliography}
\end{document}